\documentclass[11pt,a4paper]{article}
\usepackage[margin=1in]{geometry}
\usepackage{latexsym, amsfonts, amsmath}
\newcommand{\be}{\begin{equation}}
\newcommand{\ee}{\end{equation}}
\newcommand{\bea}{\begin{eqnarray}}
\newcommand{\eea}{\end{eqnarray}}
\newcommand{\bean}{\begin{eqnarray*}} 
\newcommand{\eean}{\end{eqnarray*}}

\newcommand{\brray}{\begin{array}}
\newcommand{\erray}{\end{array}}
\newcommand{\ben}{\begin{equation}{nonumber}}
\newcommand{\een}{\end{equation}{nonumber}}

\newtheorem{dfn}{Definition}[section]
\newtheorem{thm}[dfn]{Theorem}
\newtheorem{lmma}[dfn]{Lemma}
\newtheorem{ppsn}[dfn]{Proposition}
\newtheorem{crlre}[dfn]{Corollary}
\newtheorem{xmpl}[dfn]{Example}
\newtheorem{rmrk}[dfn]{Remark}
\newcommand{\bdfn}{\begin{dfn}}
\newcommand{\bthm}{\begin{thm}}
\newcommand{\blmma}{\begin{lmma}}
\newcommand{\bppsn}{\begin{ppsn}}
\newcommand{\bcrlre}{\begin{crlre}}
\newcommand{\bxmpl}{\begin{xmpl}}
\newcommand{\brmrk}{\begin{rmrk}}
\newcommand{\edfn}{\end{dfn}}
\newcommand{\ethm}{\end{thm}}
\newcommand{\elmma}{\end{lmma}}
\newcommand{\eppsn}{\end{ppsn}}
\newcommand{\ecrlre}{\end{crlre}}
\newcommand{\exmpl}{\end{xmpl}}
\newcommand{\ermrk}{\end{rmrk}}

\newcommand{\IC}{\mathbb{C}}

\newcommand{\IR}{\mathbb{R}}

\newcommand{\cla}{{\cal A}}
\newcommand{\clb}{{\cal B}}
\newcommand{\clc}{{\cal C}}
\newcommand{\cld}{{\cal D}}
\newcommand{\cle}{{\cal E}}
\newcommand{\clf}{{\cal F}}

\newcommand{\clh}{{\cal H}}
\newcommand{\cli}{{\cal I}}

\newcommand{\cll}{{\cal L}}

\newcommand{\cln}{{\cal N}}

\newcommand{\clp}{{\cal P}}
\newcommand{\clq}{{\cal Q}}

\newcommand{\cls}{{\cal S}}

\newcommand{\clu}{{\cal U}}

\def\a*{{\cal A}_{h,*}}
\def\B{{\cal B}(h)}
\def\B1{{\cal B}_1(h)}
\def\b{{\cal B}^{\rm s.a.}(h)}
\def\b1{{\cal B}^{\rm s.a.}_1(h)}

\newcommand{\ot}{\otimes}

\newcommand{\raro}{\rightarrow}

\newcommand{\lgl}{\langle}
\newcommand{\rgl}{\rangle}

\def \qed {$\Box$}

\def\a*{{\cal A}_{h,*}}
\def\B{{\cal B}(h)}
\def\B1{{\cal B}_1(h)}
\def\b{{\cal B}^{\rm s.a.}(h)}
\def\b1{{\cal B}^{\rm s.a.}_1(h)}

\begin{document}
\begin{center}
{\Large{\bf Non-existence of faithful isometric action of compact quantum groups on compact, connected Riemannian
manifolds }}\\ 
{\large {\bf Debashish Goswami\footnote{Partially supported by Swarnajayanti Fellowship and J C Bose Fellowship (2016 onwards)  from D.S.T. (Govt. of India) and also acknowledges
 the Fields Institute, Toronto for providing hospitality for a brief stay when a small part of this work was done.} and \bf Soumalya Joardar \footnote{Acknowledges support 
 from CSIR and NBHM.}}}\\ 
Indian Statistical Institute\\
203, B. T. Road, Kolkata 700108\\
Email: goswamid@isical.ac.in,\\ 
Phone: 0091 33 25753420, Fax: 0091 33 25773071\\
\end{center}
\begin{abstract}
Suppose that a compact quantum group 
$\clq$ acts faithfully  on a smooth, compact, connected  manifold
$M$, i.e. has a $C^{\ast}$ (co)-action $\alpha$ on $C(M)$, such that  the 
action $\alpha$ is  isometric in  the sense of \cite{Goswami} for some Riemannian structure on $M$.
 We prove that   $\clq$ must be commutative as a $C^{\ast}$ algebra i.e.
$\clq\cong C(G)$ for some compact group $G$ acting smoothly on $M$. In particular, the  quantum isometry group of $M$ (in the sense of \cite{Goswami}) coincides with $C(ISO(M))$. 
 
\end{abstract}
 \section{Introduction}
 It is  a very important and interesting problem in the theory of quantum groups and noncommutative geometry to study `quantum symmetries' of various classical and quantum structures. 
 Indeed, symmetries of physical systems (classical or quantum) were conventionally modeled by group actions, and after the advent of quantum groups, group symmetries were naturally generalized to  symmetries given by quantum group action. In this context, it is  natural to think of quantum automorphism or the full quantum symmetry groups of various mathematical and physical structures. The underlying basic
 principle of defining a quantum automorphism group of a given mathematical structure consists of two steps : first, to identify (if possible) the group of automorphisms of
 the structure as a universal object in a
   suitable category, and then, try to look for the universal object in a similar but bigger category by replacing groups by quantum groups of appropriate type. 
  
  The formulation and study of such quantum symmetries in terms of universal Hopf algebras date back to Manin \cite{manin_book}.
      In the analytic set-up 
    of compact quantum groups, it was considered by S. Wang who defined and studied quantum permutation groups of finite sets and quantum automorphism groups 
    of finite dimensional algebras. Subsequently, such questions were taken up by a number of mathematicians including Banica, Bichon, Collins and others
     (see, e.g. \cite{ban_1}, \cite{bichon}, \cite{wang}), 
    and more recently in the framework of Connes' noncommutative geometry (\cite{con}) by Goswami, Bhowmick,  Skalski, Banica, Soltan, De-Commer, Thibault 
    and many others who have extensively studied 
    the quantum group of isometries  (or quantum isometry group) defined in \cite{Goswami}  (see also \cite{skalski_bhow}, \cite{Soltan} etc.). 
    
    In this context, 
    it is important to compute  quantum symmetries of classical spaces. One may hope that there are many more quantum symmetries of a given classical space than classical group
     symmetries which will help one understand the space better. Indeed, it has been a remarkable discovery of S. Wang that for $n \geq 4$, a finite set of cardinality $n$ has an infinite dimensional 
      compact quantum group $\cls_n^+$ (`quantum permutation group') of symmetries. For the relevance of quantum group symmetries in a wider and more geometric context, we refer the reader to 
    the discussion on `hidden symmetry in algebraic geometry' in Chapter 13 of \cite{manin_book} where Manin made a remark about possible genuine 
    Hopf algebra symmetries of classical smooth algebraic varieties. 
   
   From a topological and geometric point of view, it is certainly more interesting to look for  examples of  faithful continuous action by genuine (i.e. not of the form 
    $C(G)$ for a compact group $G$) compact quantum groups on $C(X)$ for a connected compact  space $X$.  Several examples of such actions on 
      have been  
 constructed by H. Huang in \cite {Huchi}.  One can also adapt an  example of \cite{Etingof_walton_1} to get a faithful action of
  the group $C^*$ algebra $C^*(S_3)$ of  the group of permutations of $3$ objects (which is a finite dimensional compact quantum group)  on $C(X)$ where $X$ is the wedge product 
   of two copies of $[-1,1]$ identifying the point $0$.  
  However,  none of the above examples  are smooth manifolds. Motivated by the fact that
  a topological action of compact group on a  smooth manifold is smooth if and only if it is 
    isometric w.r.t. some Riemannian structure on the manifold, the first author of this paper and some of his collaborators and students 
     tried to compute quantum isometry groups for several classical (compact)  Riemannian manifolds including the spheres and the tori. 
   Quite remarkably,  in each of  these cases, the quantum isometry group turned out to be the same as $C(G)$ where $G$ is the corresponding isometry group.  
On the other hand,  Banica et al (\cite{banica_jyotish}) ruled out the possibility of (faithful) isometric actions of a 
 large class of compact quantum groups including $\cls_n^+$ on a connected compact Riemannian manifold. All these led the first author of the present paper to make the 
following conjecture in \cite{Rigidity}, where he also gave some supporting evidence to this conjecture considering certain class of homogeneous spaces.\\   
{\bf Conjecture I:} {\it It is not possible to have smooth faithful action of a 
genuine compact quantum group on $C(M)$ when $M$ is a compact connected smooth
manifold.}\\
The aim of this article is to prove a very important case of this conjecture, namely, under the condition that the action is isometric in the sense of \cite{Goswami} for some 
 Riemannian metric on the manifold. 
 
More precisely, we have:\\
{\bf Theorem I}\\ {\it Suppose that a CQG 
$\clq$ acts faithfully  on a smooth, compact connected  manifold
$M$ such that  the 
action is  isometric with respect to some Riemmannian structure on $M$.  
 Then  $\clq$ must be commutative as a $C^{\ast}$ algebra i.e.
$\clq\cong C(G)$ for some compact group $G$ acting smoothly on $M$.}

 
 \brmrk
 Smoothness of $M$ and compactness of the quantum group $\clq$ are quite crucial for the above conjecture. We already mentioned counter-examples in case 
  the space is non-smooth and let us refer the reader to \cite{Wal_Wang}(Example 2.20) as well as \cite{loc_comp}
  for faithful algebraic (co)-action of Hopf-algebras corresponding to genuine non-compact 
  quantum groups on commutative domains associated with affine varieties.
  However,  we do not yet have any 
example of genuine CQG action on non-compact, smooth, connected manifold and believe 
    that it may be possible to extend our no-go result to this case as well. 
 \ermrk

  \brmrk  
    In some sense, our results indicate that one cannot possibly have a genuine `hidden quantum symmetry' in the sense of Manin (Chapter 13 of \cite{manin_book}) 
     for smooth connected varieties coming from CQG Hopf algebras; i.e. one must look for such quantum symmetries given by Hopf algebras of non-compact type only.
      From a physical point of view, it follows that for a classical mechanical  system with phase-space modeled by a compact connected manifold,
 the generalized notion of symmetries in terms of (compact) quantum groups coincides with the conventional notion, i.e. symmetries coming from group actions. \ermrk

\brmrk In \cite{Etingof_walton_1}, Etingof and Walton  obtained
a somewhat similar result in the purely algebraic set up of finite dimensional
Hopf algebra actions on commutative domains. However, their proof does not seem to extend to the infinite dimension as it crucially
depends on semisimplicity and finite dimensionality of the Hopf algebra.  \ermrk
\brmrk
  We should also mention here the  attempts by several authors  to formulate  a notion of quantum isometry group in the purely metric space set-up 
  (see \cite{sabbe}, \cite{ban_1}, \cite{Metric} etc.) and the proof  by A. L. Chirvasitu (\cite{chirvasitu}) of 
  non-existence of genuine quantum isometry in the metric space set-up for geodesic metric of negatively curved, compact connected Riemannian manifold.
  \ermrk

 The following fact, which was  observed in earlier works like
\cite{Rigidity}, plays a crucial role in the proof of Theorem I.\\
{\it Fact :} There does not exist  any faithful action   by a genuine CQG on a subset with
nonempty interior of some Euclidean space $\IR^n$ which is affine in the sense
that the action leaves  the linear span of the coordinate functions and
the constant function $1$ invariant.\\
 We prove  Theorem I by deducing first  that the given isometric CQG action can be  lifted  to a faithful affine action on the closure of a suitable 
  bounded Euclidean domain (open connected set with smooth boundary).
 In the classical case, i.e. a compact group action, one may take bases in sufficiently many 
 spectral subspaces of $C^\infty(M)$ and use functions in these bases to embed $M$ into $\IR^N$ for some $N$, so that the action becomes affine (equivalently linear, after a suitable shift 
  of the 
  origin) w.r.t. the coordinates of $\IR^N$. This can be done for CQG actions as well. More precisely, the action (say $\alpha$) will satisfy 
   \be \label{bbb} \alpha(X_i)=\sum_{j=1}^N X_j \ot q_{ij} +c_i1\ee for some generating set of self-adjoint elements $q_{ij} $ of the CQG, where $c_i \in \IR$, $X_i$ is the restriction of the 
    $i$-th coordinate function of $\IR^N$ to the subset $M \subset \IR^N$. However,  the main difference 
     between the classical and quantum situation is that an  affine representation of a group $G$ on $\IR^N$ can be lifted to an action on the algebra of continuous or smooth  
      functions on $\IR^N$ by `dualizing' the point-wise action. Using compactness of the group, we can  actually get action on some suitable  Euclidean closed ball  $B$ of large enough 
       radius $r>0$ 
       (say) 
       around $(c_1, \ldots, c_N)$ containing $M$ in the interior. This is not the case for a general CQG action. An affine representation on $\IR^N$ like the one
       obtained by restricting $\alpha$ to ${\rm Sp}\{ 1, X_1, \ldots, X_N \}$ need not induce any action on  $C(B)$. Indeed, existence of   such an action would imply  in particular 
        that 
         the algebra generated by 
        $\{ \alpha(\tilde{X}_i)(p), p\in B,~~i=1, \ldots, N\}$ is commutative, where $\tilde{X}_i$ is the restriction of the $i$-th coordinate function to $B$. 
        This is equivalent to the commutativity of the 
        algebra generated by $\{ \sum_j p_jq_{ij}:~p\equiv (p_1, \ldots, p_N)\in B,~i=1, \ldots, N \} $. 
        However, the fact that $\alpha$ is a $\ast$-homomorphism on $C^\infty(M)$  gives us only 
         the commutativity of $\{ \sum_j p_jq_{ij}:~(p_1, \ldots, p_N)\in M,~i=1, \ldots, N \} $. Moreover, the coordinate functions restricted to $M$ need not be `quadratically independent' in the sense
          of \cite{Rigidity}, so that we cannot  deduce the commutativity of $q_{ij}$'s. 
          
   For this reason, we need to take a more elaborate route for lifting the given CQG action to an affine action on a bounded domain.  

       (a) We begin by  lifting the isometric action to the total space $O(M)$ of the orthonormal frame bundle.\\
        (b) As $O(M)$ is parallelizable, one can choose an isometric embedding of it in some Euclidean space 
      where the corresponding normal bundle is trivial. Thus, the isometric action  on $O(M)$ can be further lifted to 
       an isometric action on some suitable tubular neighbourhood, say $N$, of $O(M)$.\\ 
      (c) Finally, as  $N$ is locally isometric to the flat Euclidean space, any isometric action on it must be affine 
      in the corresponding coordinates. The connectedness of $M$ is used only in this step.
 
 In some sense, the liftings in the above steps are achieved by adapting the classical line of  arguments to the CQG set-up. 
  However, this adaptation is  rather non-trivial due to 
noncommutativity at every stage. Indeed, to start with, we only have the commutativity of $\alpha(f)(m)$ with $\alpha(g)(m)$ for different $f,g \in C^\infty(M)$, where $m \in M$ is fixed.
 Using the isometry condition, we first deduce that the first order partial derivatives of $\alpha(f)$ at $m$ (with respect to 
  any local coordinates) will commute among themselves as well as with $\alpha(g)(m)$ for 
  $f,g \in C^\infty(M)$. This allows us to prove that the natural analogue of the differential of the action $\alpha$, which is a representation of the CQG on the module of one and higher 
   forms, 
   is well-defined. This is necessary to lift the action to $O(M)$. Then we prove that in a suitable sense the reprsentaion on the module of one-forms 
   `commutes' with the Levi-civita connection for  the Riemannian metric. This in turn helps us to deduce `higher order commutativity', namely the commutativity 
    of all order partial derivatives of $\alpha(f)(m)$ $f \in C^\infty(M)$  at a given point $m$ of the manifold. This fact has been crucially used in the proof of the step (c).

 We begin by collecting some basic definitions, notations and standard facts  in Section
2.  In Section 3 we introduce smooth and inner product preserving action and prove that such an action can be lifted to the orthonormal frame bundle. 
 Then we discuss the implications of isometric actions (i.e. actions commuting with the Hodge Laplacian) in Section 4. 
  Finally, in
Section 5, we state and prove the main result  using embedding of $O(M)$ in some $\IR^N$ with trivial normal bundle and the results of previous sections.

\section{Preliminaries}
\subsection{Notational convention} In this paper all the
Hilbert spaces are over $\mathbb{C}$ unless mentioned otherwise. If $V$ is a
vector space over real numbers we denote its complexification by
$V_{\mathbb{C}}$. For a vector space $V$, $V^{\prime}$ stands for its algebraic dual. We denote the domain of a linear (possibly unbounded) map $L$ by ${\rm Dom}(L)$. 
The algebraic tensor product of modules over an algebra $\clc$ is denoted by $\ot_\clc$. In $\clc=\IC$, i.e. the modules are vector spaces, we simply use $\ot$. 
 For a complex $\ast$-algebra $\clc$, let $\clc_{s.a.}=\{c\in\clc: c^{\ast}=c\}$. We shall denote the $C^{\ast}$
algebra of bounded operators on a Hilbert space $\clh$ by $\clb(\clh)$ and the
$C^{\ast}$ algebra
of compact operators on $\clh$ by $\clb_{0}(\clh)$. $Sp$, $\overline{Sp}$ 
stand for the
linear span and closed linear span of elements of a vector space respectively, whereas ${\rm Im}(A)$ denotes the image of a linear map. $\delta_{ij}$ is the Kronecker delta function.
For an algebra  $\clc$, 
 let 
$\sigma_{ij}:\underbrace{\clc\ot\clc\ot...\ot\clc}
_{n-times}\raro\underbrace{\clc\ot\clc\ot...\ot\clc}_{n-times}$ (with $\sigma_{12} \equiv \sigma $ for $n=2$) denote  the  map which flips 
the  $i$ and $j$-th tensor copies.

We'll use the notation $\hat{\otimes}$ for the minimal (spatial) tensor product between $C^*$ algebras or tensor product between $\ast$-homomorphisms (more generally, 
 completely positive maps) between $C^*$ algebras.  We will usually denote scalar valued inner product of Hilbert spaces by $<\cdot, \cdot >$ and some $\ast$-algebra 
valued inner product of Hilbert modules by $\lgl\lgl \cdot, \cdot \rgl\rgl$.
For two  Hilbert $\cla$-$\clb$  and $\clb$-$\clc$ bimodules $\cle$ and $\cle^\prime$, (hence in particular for Hilbert spaces) 
  we denote 
  by $\cle \overline{\otimes} \cle^\prime$ the Hilbert $\cla$-$\clc$ bimodule obtained by suitably quotienting and completing 
    the algebraic tensor product $\cle \ot_\clb \cle^\prime$
   w.r.t. the $\clc$-valued inner product defined by  
   $\lgl\lgl \xi \ot_\clb \xi^\prime, \eta \ot_\clb \eta^\prime\rgl\rgl=\lgl\lgl\xi^\prime, \lgl\lgl \xi, \eta\rgl\rgl\eta^\prime \rgl\rgl$ on the simple tensor elements
    and then extended naturally to their linear span. 
   We'll have occasions to use Hilbert bimodules over certain 
    Fre\'chet $\ast$-algebras too.

  \subsection{Compact quantum groups, their representations and actions}
We very briefly outline the notion of compact quantum groups (CQG) and their representations. The reader is referred to \cite{Van}, \cite{Pseudogroup} and the references therein for details.
A compact quantum group (CQG for short) is a  unital $C^{\ast}$-algebra $\clq$ with a
coassociative coproduct 
(see \cite{Van}) $\Delta$ from $\clq$ to $\clq \hat{\ot} \clq$  
such that each of the linear spans of $\Delta(\clq)(\clq\ot 1)$ and that
of $\Delta(\clq)(1\ot \clq)$ is norm-dense in $\clq \hat{\ot} \clq$. 
From this condition, one can obtain a canonical dense unital $\ast$-subalgebra
$\clq_0$ of $\clq$ on which linear maps $\kappa$ and 
$\epsilon$ (called the antipode and the counit respectively) are defined making the above subalgebra a Hopf $\ast$-algebra.

 It is known that there is a unique state $h$ on a CQG $\clq$ (called the Haar
state) which is bi invariant in the sense that $({\rm id} \ot h)\circ
\Delta(a)=(h \ot {\rm id}) \circ \Delta(a)=h(a)1$ for all $a\in\clq$. The Haar state
need not be faithful in general, though it is always faithful on $\clq_0$ at
least. The image of $\clq$ in the GNS representation of $h$, say $\clq_{r}$, is called the reduced CQG corresponding to $\clq$. 
 
 A unitary representation of a CQG $(\clq,\Delta)$ on a Hilbert
space $\clh$ is a $\mathbb{C}$-linear map $U$ from $\clh$ to the Hilbert module
$\clh\bar{\ot}
\clq$ such that \\
1. $\lgl\lgl U(\xi),U(\eta)\rgl\rgl=\lgl\xi,\eta\rgl 1_{\clq}$, where $\xi,\eta\in \clh$.\\
2. $(U\overline{\ot} {\rm id})U=({\rm id}\hat{\ot} \Delta)U,$\\
3. The right $\clq$-linear span of ${\rm Im}(U)$ is dense in $\clh \overline{\ot} \clq$.\\
 \indent In 2. above, the map $(U \overline{\ot} {\rm id})$ denotes the extension of $U \ot {\rm id}$ to the completed tensor product 
  $\clh \overline{\ot} \clq$ which exists by the isometry condition 1. A closed subspace $\clh_{1}$ of $\clh$ is said to be invariant if
$U(\clh_{1})\subset \clh_{1}\bar{\ot}\clq$. 

\bdfn
An algebraic (co)-representation (or co-module) for a Hopf algebra $(H,\Delta,\epsilon,\kappa)$ on  a vector space $V$ is  
a $\IC$-linear map  $\alpha : V \raro V \ot H$ such that $(\alpha \ot {\rm id}) \circ \alpha=({\rm id} \ot \Delta) \circ \alpha$ and $({\rm id} \ot \epsilon) \circ \alpha={\rm id}$.
In case $H$ is  a Hopf $\ast$ algebra, and $V=A$ is a unital $\ast$ algebra and the (co)-representation $\alpha$ is also a unital $\ast$ algebra homomorphism,
 we say that it is a  Hopf $\ast$ algebraic (co)-action of $H$ on $A$. 
\edfn
For an algebraic (co)-representation  as above, one can prove that ${\rm Sp} ( \alpha(V) H)=V\ot H$. For a Hopf algebra $H$ with the coproduct $\Delta$, 
 we write $\Delta(q)=q_{(1)}\ot q_{(2)}$ suppressing the summation notation (Sweedler's notation). For
an algebra (other than $H$ itself) or module $\cla$ and a $\mathbb{C}$-linear map $\Gamma:\cla\raro
\cla\ot H$ (typically a comodule map or a coaction) we shall also use an analogue of Sweedler's notation and write
$\Gamma(a)=a_{(0)}\ot a_{(1)}$. 

\bdfn
\label{CQG_action}
A unital $\ast$-homomorphism  $\alpha:\clc\raro \clc\hat{\ot}\clq$, where $\clc$ is a unital $C^\ast$-algebra and 
 $\clq$ is a CQG, is
said to be an action of $\clq$ on $\clc$  if\\
1. $(\alpha \hat{\ot} {\rm id})\alpha=({\rm id} \hat{\ot} \Delta)\alpha$ (co-associativity).\\
2. $Sp \ \alpha(\clc)(1\ot \clq)$ is norm-dense in $\clc\hat{\ot}\clq$. 

We say that the action is 
algebraic over a $\ast$-subalgebra $\clc_{0}\subset\clc$ if $\alpha|_{\clc_{0}}:\clc_{0}\raro\clc_{0}\ot\clq_{0}$ is a Hopf $\ast$-algebraic (co)-action.
\edfn
The following result is proved by Podles (\cite{Podles}).
\bppsn
\label{maximal}
 Given an action $\alpha$ of a CQG $\clq$ on $\clc$, there exists a norm-dense unital $\ast$-subalgebra  $\clc_{0}$ of $\clc$ 
 over which $\alpha$ is algebraic.
\eppsn
The action $\alpha$ is said to be faithful if the $\ast$-subalgebra
of $\clq$ generated by the elements of the form $(\omega\ot {\rm id})(\alpha(a))$, $a \in \clc_0, ~\omega \in \clc_0^\prime$
 is norm-dense in $\clq$. 
When $\clc=C(M)$, faithfulness is equivalent to requiring the density of the subalgebra generated by elements of the form $\alpha(f)(m)$, $f \in C(M)$, $m \in M$.
\subsection{$C^{\infty}(M)$ and bimodule of forms}
\label{form}
Let $M$ be a smooth, $n$-dimensional compact manifold possibly with boundary. We denote the algebra of real (complex respectively) 
valued smooth functions on $M$ by $C^{\infty}(M)_{\mathbb{R}}$ ($C^{\infty}(M)$ respectively). The natural Fr\'echet topology on $C^{\infty}(M)$ is 
given by the seminorms of the form $p^{U, K,\alpha}$,
\begin{displaymath}
p^{U,K,\alpha}(f)={\rm sup}_{x\in K}|\partial^{\alpha}(f)(x)|,
\end{displaymath}
where $K$ is a compact subset contained in the domain of some coordinate chart $(U,(x_{1},...,x_{n}))$, $\alpha=(i_{1},...,i_{k})$ a multi index and 
$\partial^{\alpha}=\frac{\partial}{\partial x_{i_{1}}}...\frac{\partial}{\partial x_{i_{k}}}$, $i_j \in \{1, \ldots, n\}.$.
We can similarly define a Fr\'echet topology on $C^{\infty}(M,E)$, the space of smooth $E$-valued functions on $M$ for any Fr\'echet space $E$. 
A word on our notational convention: we denote by $T_{m}^{\ast}(M)$ the complexified cotangent space at $m$, whereas $(T^{\ast}_{m} M)_{\IR}$ will denote the corresponding 
real vector space. 
 Let
$\Omega^{1}(C^\infty(M)) \equiv \Lambda^{1}(C^{\infty}(M))$ be the space of smooth complex valued $1$ forms on the manifold
$M$, with the natural locally convex
topology induced by the topology of $C^{\infty}(M)$ given by a
family of seminorms $q^{U,K,\alpha}(\omega)={\rm sup_{x\in K, 1\leq j\leq n}}|\partial^{\alpha}f_{j}(x)|$, 
where $K \subset U$, $\alpha$ are as before, 
and $\omega|_{U}=\sum_{j=1}^{n}f_{j}dx_{j}$. It is clear from the definition that the differential map
$d:C^{\infty}(M)\raro \Omega^{1}(C^{\infty}(M))$ is Fr\'echet continuous. As $M$ is compact, there is a Riemannian structure. Using the Riemannian
structure on $M$ we can equip $\Omega^{1}(C^{\infty}(M))$ with a $C^{\infty}(M)$
valued inner product given by $\lgl\lgl\omega,\eta\rgl\rgl(m)=\lgl\omega(m),\eta(m)\rgl_{m}$, where $\lgl,\rgl_{m}$ denotes the complex-valued  inner product on complexified cotangent 
space coming from the Riemannian structure. This makes  $\Omega^{1}(C^{\infty}(M))$  a Fr\'echet-Hilbert module over $C^\infty(M)$. 

It can be shown  that $\Omega^{1}(C^{\infty}(M))= \ {\rm Sp} \{fdg:f,g\in C^{\infty}(M)\}$.

Given $\omega \in \Omega^1(M)$, let us denote by $X_\omega$ the (unique) smooth vector field satisfying
\begin{eqnarray}
\label{x_omega}
 X_\omega(f)=\lgl\lgl \omega, df\rgl\rgl
\end{eqnarray}
 for any smooth function $f$.
 Moreover, one has $\lgl\lgl X_\omega, X_\eta\rgl\rgl=\lgl\lgl\eta, \omega \rgl\rgl$ from the relation between the Riemannian 
  inner product on the tangent and cotangent spaces.

Similarly, we construct  the Fr\'echet-Hilbert bimodules $\Omega^k(C^\infty(M))=\Omega^1(C^\infty(M)) \otimes_{C^\infty(M)}  \ldots \otimes_{C^\infty(M)} \Omega^1(C^\infty(M))$ 
by taking $k$-fold tensor product over $C^\infty(M)$ and the inner product (denoted by $\lgl\lgl \cdot, \cdot \rgl\rgl$) coming from the Riemannian structure. In fact, we have $\lgl\lgl\omega_{1}\ot_{C^{\infty}(M)}...\ot_{C^{\infty}(M)}\omega_{k},\eta_{1}\ot_{C^{\infty}(M)}...\ot_{C^{\infty}(M)}\eta_{k}\rgl\rgl=\lgl\lgl\omega_{1},\eta_{1}\rgl\rgl...\lgl\lgl\omega_{k},\eta_{k}\rgl\rgl$, using the commutativity of $C^{\infty}(M)$ and the fact that $\omega f=f\omega$ for all $\omega\in \Omega^{1}(C^{\infty}(M))$ and $f\in C^{\infty}(M)$.
Clearly, the group $S_k$ of permutations of $k$ objects acts on $\Omega^k(C^\infty(M))$ in an obvious way by 
 permuting the tensor copies. For a character $\chi$ of $S_k$, let $P_\chi$ be the corresponding spectral projection given by:
 $$ P_{\chi}(f \omega_1 \ot \omega_2 \ot \ldots \ot \omega_k)=\frac{1}{k!} \sum_{ \sigma \in S_k} \chi(\sigma) f \omega_{\sigma(1)} \ot \omega_{\sigma(2)} \ot \ldots 
   \ot \omega_{\sigma(k)},$$
  where $f \in C^\infty(M), \omega_i \in \Omega^1(C^\infty(M))~\forall i$. 
 The bimodule of  smooth 
$k$ forms $\Lambda^{k}(C^{\infty}(M))= \ {\rm Sp} \ \{f df_{1}\wedge...\wedge df_{k}:f, f_{i}\in C^{\infty}(M)\}$ is a complemented  submodule $P_{\rm sgn}(\Omega^k(C^\infty(M))$ 
 of $\Omega^k(C^\infty(M))$, comprising the antisymmetric part, where ${\rm sgn}( \sigma)$ denotes the sign of the permutation $\sigma$. 
  It is easy to see that the $S_k$-action preserves the inner product on the product Hilbert space. 
  Clearly,  the above arguments go through if we 
replace $C^{\infty}(M)$ by any unital $\ast$-subalgebra $\cla$.
In fact, the $S_k$-action comes from the fibre-wise permutation on the bundle $T^{\ast}M \ot T^{\ast}M \ot \ldots \ot T^{\ast}M$ ($k$ copies). It also preserves the 
 inner product of  the product Hilbert space structure coming from the inner product on $T^{\ast}_{m}M$ given by the Riemannian metric, hence becomes a unitary representation on each fibre. 
  From this, it follows easily that the idempotent $P_\chi$'s are orthogonal projections satisfying $\sum_\chi P_\chi ={\rm Id}$, $P_\chi P_{\chi^\prime}=0$ if $\chi \ne \chi^\prime$. 
  Therefore, we have a decomposition of $\Omega^k(\cla)$  into mutually orthogonal (w.r.t.  $\lgl\lgl \cdot, \cdot\rgl\rgl$) 
  subspaces $P_\chi(\Omega^k(\cla))$.
  In the  special case $k=2$, we write $\clf_s(\cla)= (I-P_{\rm sgn})(\Omega^2(\cla))$ (`symmetric submodule')
   and $\Lambda^2(\cla)=P_{\rm sgn}(\Omega^2(\cla))$ (`antisymmetric submodule'). We have  $\Omega^2(\cla)=\clf_s(\cla)\oplus \Lambda^2(\cla)$, 
    so that $\Lambda^2(\cla)=\clf_s(\cla)^\perp:=\{ \omega \in \Omega^2(\cla):~\lgl\lgl\omega, \eta\rgl\rgl=0~\forall \eta \in \clf_s(\cla) \}$. 
    We also note the following equivalent description of $\clf_s(\cla)$ (see, e.g. \cite{Landi}, page 101):
\bppsn 
\label{symm_description}
The submodule $\clf_s$  
    is spanned (as a right $\cla$-module) by $\sum_i (df_i \ot_\cla dg_i)$, with $f_i,g_i$ such that $\sum_i f_idg_i=0$ (finite sum).\eppsn

We remark that if $\cla$ is a Fr\'echet dense subalgebra of $C^{\infty}(M)$, then $\Lambda^{k}(\cla):={\rm Sp}\{f df_{1}\wedge...\wedge df_{k}:f, f_{i}\in\cla\}$ is 
dense in the Fr\'echet Hilbert module $\Lambda^{k}(C^{\infty}(M))$ for all $k=1,...,n.$\\
Let $E$ be a smooth, hermitian vector bundle over $M$ and $\cle$ be the $C^\infty(M)$-module of smooth sections of $E$. 
 For a $C^{\ast}$ algebra $\clq$, there is a natural $C^\infty(M, \clq)$-bimodule structure as well as $C^\infty(M, \clq)$-valued inner product $\lgl\lgl \cdot, \cdot \rgl\rgl$ 
 on $\cle \ot \clq$, satisfying 
  the following:
	\begin{displaymath}
 ( f\ot q_1)(s \ot q^\prime)(g \ot q_2)=fsg \ot q_1q^\prime q_2,~~\lgl\lgl s \ot q, s^\prime \ot q^\prime\rgl\rgl=\lgl\lgl s,s^\prime\rgl\rgl_\cle \ot q^*q^\prime,
\end{displaymath} where 
 $s,s^\prime \in \cle, $ $f,g \in C^\infty(M),$ $q,q^\prime, q_1,q_2 \in \clq$ and $\lgl\lgl \cdot, \cdot\rgl\rgl_\cle$ denotes the $C^\infty(M)$-valued inner product on $\cle$ coming 
  from the hermitian structure of $E$.
 We topologize $\cle \ot \clq$ with the weakest locally convex topology making the map $\Psi \mapsto \lgl\lgl\Psi, \Psi \rgl\rgl \in C^\infty(M, \clq)$ continuous with respect to 
  the Fr\'echet topology. The corresponding completion of $\cle \ot \clq$ will be denoted by $\cle \bar{\ot} \clq$ and it  also inherits by continuity  the $C^\infty(M,\clq)$ bimodule 
  structure as well as the $C^\infty(M, \clq)$-valued inner product. 
  
  An element $X \in \cle \bar{\ot} \clq$ can be thought of as a smooth, $\clq$-valued section. To see this, fix $m \in M$ and choose smooth sections $s_1, \ldots, s_k$ (where $k$ is the 
   rank of $E$) such that $\{s_1(m), \ldots, s_k(m)\}$ is an orthonormal basis of the fibre $E_m$ of $E$ at $m$. We define $X(m):=\sum_i s_i(m) \ot X_i(m)  \in E_m \ot \clq$,
    where $X_i(m)=\lgl\lgl s_i \ot 1, X\rgl\rgl(m)$. It is easily seen that  the definition does not depend on the choice of $s_{1},...,s_{k}$, except their values at $m$. One can check this 
      in case $X$ is chosen from the algebraic tensor product $\cle \ot \clq$, and by density, it follows for a general $X$. In fact, we have the following:
      \blmma
     There is a one-to-one  correspondence between elements of $\cle \overline{\ot} \clq$ and the smooth $\clq$-valued sections, i.e. maps $X : M \raro \bigcup_m E_m \ot \clq$,
      such that $X(m) \in E_m \ot \clq$ for all $m$, and for every $\xi \in \cle$, $m \mapsto \lgl\lgl \xi(m) \ot 1, X(m) \rgl\rgl \in \clq$ is smooth. 
      \elmma
   {\it Proof:}\\ We have already seen one direction of this statement. To see the reverse direction, let $X$ be a smooth $\clq$-valued section of the bundle $E$. Choose 
   a finite open cover $U_{\alpha}, \alpha =1,\ldots, l$ (say) such that $E|_{U_\alpha}$ is trivial. Let $f_{\alpha},\alpha=1, \ldots, l$ be the corresponding smooth partition of unity and 
   for each $\alpha$, choose smooth sections   $\{ s_i^\alpha, i=1, \ldots, k \}$ such that $\{ s_i^\alpha(m) , i=1, \ldots ,k\}$ is a basis of $E_m$ $\forall m \in U_\alpha$.
   We can write $X(m)$ in terms of this basis, say, $X(m)=\sum_i s^{\alpha}_{i} (m) \ot q^{\alpha}_{i}(m)$, where $q^{\alpha}_{i}(m) \in \clq$. 
   It is easily seen that $m \mapsto q^{\alpha}_{i}(m)$ defines an element of $C^{\infty}(M, \clq)$. Hence, $X=\sum_{\alpha, i} f_\alpha (s^{\alpha}_{i} \ot 1)q^{\alpha}_{i}$, which is an element 
    of $\cle \bar{\ot} \clq$. \qed

\indent In particular,  we will 
identify elements $\Omega\in\Omega^{1}(C^{\infty}(M))\bar{\ot}\clq$ with $\clq$-valued smooth one forms, i.e. $\Omega:M\raro \cup_{m\in M}(T^{\ast}_{m}M)\ot\clq$, 
such that  $\Omega(m)\in T^{\ast}_{m}M\ot\clq$ for all $m\in M$ and for any local coordinate chart  $(U,(x_{1},...,x_{n}))$ around $m\in M$, 
we have  $\Omega(x)=\sum_{i=1}^{n}dx_{i}(x)\ot\Omega_{i}(x)~\forall x \in U$,  for some  $\Omega_{i}\in C^{\infty}(M,\clq)$, $i=1,...,n.$
We will usually write $dx_{i}(x)\ot\Omega_{i}(x)$ as $dx_{i}(x)\Omega_{i}(x)$ and $\Omega=\sum_{i=1}^{n}dx_{i}\Omega_{i}$ on $U$. 

For  a smooth vector field $X$ defined on some open subset ${\rm Dom}(X)$  of $M$  and $F \in C^\infty(M, \clq)$, 
 we define $X(F)(m):=\frac{d}{dt}|_{t=0} F(\gamma(t))$, where $m \in {\rm Dom}(X)$ and  $\gamma$ is the integral curve for $M$ passing through $m$. We also  define $dF\equiv (d\ot {\rm id})(F)\in\Omega^{1}(C^{\infty}(M))\bar{\ot}\clq$ for $F\in C^{\infty}(M,\clq)$ by the following:
\begin{displaymath}
 (dF)(m):= \sum_{i=1}^{n} dx_{i}(m)(\frac{\partial F}{\partial x_{i}})(m),
\end{displaymath}
for $m\in M$ and for any local coordinate chart $(U,(x_{1},...,x_{n}))$ around $m$. Clearly this is uniquely defined by the condition 
\begin{displaymath}
 \omega(dF(m))=d(F_\omega)(m)
\end{displaymath}
for every bounded linear functional $\omega$ on $\clq$,  where $F_\omega \in C^\infty(M)$ is given by $F_\omega(m):=\omega(F(m))$. Thus $dF$ does not depend on the choice of the local coordinates.

Similarly, we consider $\clq$-valued $k$-forms $\Lambda^k(M, \clq) \equiv \Lambda^k(C^\infty(M)) \overline{\otimes} \clq$. 
     We also  define  the differential and the wedge product of $\clq$-valued forms. For $\Omega \in \Lambda^k(M, \clq), \Theta \in \Lambda^l(M, \clq)$, we have 
$\Omega \wedge \Theta \in \Lambda^{k+l}(M, \clq)$ given by $(\Omega \wedge \Theta)(m):=\Omega(m) \wedge \Theta(m)$, 
Where $\wedge : (\IC^k \otimes \clq ) \times (\IC^l \ot \clq) \raro (\IC^k \wedge \IC^l) \ot \clq$  
is defined by $(v \ot q) \wedge (v^\prime \ot q^\prime):=v \wedge v^\prime \ot qq^\prime.$
 
 For a $k$-form $\omega$, we denote by $\overline{\omega}$ the form obtained by fibre-wise complex conjugation. Let $\Lambda^k(M)_\IR$
  denote the module of forms $\omega$ such that $\overline{\omega}=\omega$. This is a module over $C^\infty(M)_\IR$.  We also have an analogue of complex conjugation 
   on $\Lambda^{k}(M, \clq)$ coming from the natural conjugation on $T^{\ast}_{m}(M) \ot \clq$, namely $\overline{(v \ot q)}=\overline{v} \ot q^*$.

\subsection{\bf Basics of the normal bundle}
Let $M\subseteq \mathbb{R}^{N}$ be a smooth embedded submanifold of $\mathbb{R}^{N}$ without boundary. For each point $x\in M$ define the space of normals to $M$ at $x$ to be \begin{displaymath} N_{x}(M)=\{v\in \mathbb{R}^{N}: v \perp T_{x}(M)\}.\end{displaymath}
The total space $\cln(M)$ of the normal bundle is defined to be
\begin{displaymath}\cln(M)=\{(x,v)\in M\times \mathbb{R}^{N}; v\perp T_{x}(M)\}\end{displaymath} with the
projection $\pi$ on the first coordinate. Define
$\cln_{\epsilon}(M)=\{(x,v)\in \cln(M);||v||\leq \epsilon\}.$
It can be shown that $\cln_\epsilon(M)$ is a compact manifold of dimension $N$ with smooth boundary
$\partial \cln_\epsilon(M) =\{ (x,v) \in \cln(M):~\| v \| =\epsilon\}$ (see page no. 153
of \cite{Shastri}).
\blmma
\label{global_diff} 
(i) If a compact $n$-manifold $M$ without boundary
embedded  $\mathbb {R}^{N}$ has trivial normal bundle, then 
there exist an $\epsilon> 0$ and a global diffeomorphism $F: M\times
B_{\epsilon}^{N-n}(0)\raro \cln_{\epsilon}(M)\subseteq \mathbb{R}^{N}$   given by
$$F(x,u_{1},u_{2},...,u_{N-n})=x+ \sum_{i=1}^{N-n}\xi_{i}(x)u_{i}$$ where $B_{\epsilon}^{N-n}(0)$ is the closed  ball in $\IR^{N-n}$ of radius
$\epsilon$ centered at $0$, 
$(\xi_{1}(x),...,\xi_{N-n}(x))$ is an orthonormal basis of $N_{x}(M)$ for all
$x$, and $x\mapsto \xi_{i}(x)$ is smooth $\forall \  i=1,...,(N-n).$\\
(ii) The map 
$\pi_{F}:C^{\infty}(\cln_{\epsilon}(M))\raro
C^{\infty}(M\times B_{\epsilon}^{N-n}(0))$ given by
$\pi_{F}(f)(x,u_{1},u_{2},...,u_{N-n})=f(x+ \sum_{i=1}^{N-n}\xi_{i}(x)u_{i})$ is an algebra isomorphism. 
\elmma
{\it Proof}:\\
(i) is a consequence of the tubular neighbourhood lemma. For the proof see
\cite{Shastri}. Proof of (ii) is straightforward. \qed.
\vspace{0.2in}\\
We now introduce the notion of stably parallelizable manifolds.
\bdfn
A manifold $M$ is said to be stably parallelizable if its tangent bundle is stably trivial.
\edfn
 From the discussion following Theorem 7.2 of Chapter 9 in
\cite{Kol}  (see also  \cite{Singhof}), we get the following:
\bppsn
\label{iso_emb}
A manifold $M$ is stably parallelizable if and only if there exists an 
embedding of $M$ into some euclidean space with  trivial normal
bundle. Moreover, given a Riemannian structure on a stably parallelizable
manifold $M$, we can choose the embedding to be isometric.
\eppsn
We note that parallelizable manifolds (i.e. which have trivial tangent bundles)
are in particular stably parallelizable. Moreover, given any compact Riemannian
manifold $M$, its orthonormal frame bundle $O_{M}$ is parallelizable.

\section{Smooth and  inner product preserving  actions of a CQG on a manifold}
Throughout this section, let $M$ be a compact, smooth manifold possibly with boundary, unless otherwise mentioned.
\subsection{Smooth and inner product preserving actions}
\bdfn
\label{smooth_action}
A $C^{\ast}$-action $\alpha$ of a CQG $\clq$ on $C(M)$ is called smooth if\\ 
 (i) $ \alpha(C^\infty(M)) \subseteq C^\infty(M, \clq)$, \\
 (ii) the $\IC$-linear 
 span of $\alpha(C^\infty(M))(1 \ot \clq)$ is Fr\'echet-dense in $C^\infty(M, \clq)$.\\
 In case $M$ has a smooth boundary $\partial M$, we also require that $\alpha$ 
  maps the $C^*$ ideal $\cli:=\{ f \in C(M): f|{\partial M}=0 \}$ into $\cli \hat{\ot} \clq$. We'll often say that $\alpha$ is a smooth action of $\clq$ on $M$. 
 \edfn
 \brmrk
 It follows from the Closed Graph Theorem that any smooth action is automatically continuous with respect to the Fr\'echet topologies on $C^\infty(M)$ and $C^\infty(M, \clq)$.
 \ermrk
 An almost verbatim adaptation of arguments in  \cite{Podles} gives us the following analogue of Proposition (\ref{maximal}):
\bppsn
\label{smooth_action_dense}
A $C^{\ast}$ action $\alpha$ on $C(M)$ is smooth iff $\alpha(C^{\infty}(M))\subset C^{\infty}(M,\clq)$ and there is a Fr\'echet dense subalgebra $\cla$ of $C^{\infty}(M)$
over
which $\alpha$ is algebraic.
\eppsn
Given a smooth action $\alpha$ let us 
 introduce the following\\
{\it Notation:} For $x \in M$ let us denote by $\clq_x$ the unital $\ast$-subalgebra of
$\clq$ generated by  \be  \{ \alpha(f)(x), ((\phi \ot {\rm
id})\alpha(g))(x),~f,g \in C^\infty(M),~ \phi \in \chi(M)\}, \label{def_q_x}\ee where $\chi(M)$ is the set of smooth vector fields on $M$. Given a local coordinate 
$(x_1, \ldots, x_n)$ around a point $m$, by choosing a smooth vector field which agrees with $\frac{\partial}{\partial x_i}$ in a neighbourhood of $m$ , we see that 
$\frac{\partial}{\partial x_i}|_m \alpha(f)$ belongs to $\clq_m$.

\bdfn 
Suppose that $M$ has a Riemannian structure with the corresponding  $C^{\infty}(M)$ valued inner product $\lgl\lgl\cdot, \cdot\rgl\rgl$ on $\Omega^1(C^\infty(M))$. 
A smooth action $\alpha$ on $M$ is said to preserve  the inner product   if 
\begin{eqnarray}
\label{inner_prod_pres_111}
 \lgl \lgl d\alpha (f), d\alpha (g)\rgl\rgl=\alpha(\lgl\lgl df,dg\rgl\rgl )
 \end{eqnarray}
 for all $f,g\in C^\infty(M)$. 
\edfn
 It is easy to see, by the  Fr\'echet continuity of the maps involved that it is enough to have  (\ref{inner_prod_pres_111}) for  $f,g$ varying in some 
Fr\'echet dense unital $\ast$-subalgebra of $C^{\infty}(M)$. 
\bthm
\label{automatic_lift}

Let $\clq$ be a reduced CQG, i.e. the Haar state is faithful on $\clq$. Given a smooth action $\alpha$ of a CQG $\clq$ on $M$ the following are equivalent:\\
 (i) For every $x \in M$, we have \be \label{dalpha_eq} \alpha(f)(x) (d \ot {\rm id})(\alpha(g))(x)=(d \ot {\rm id})(\alpha(g))(x)\alpha(f)(x), \ee
 for all $f,g \in C^\infty(M)$ and  $x \in M$.\\
(ii) The manifold $M$ has a Riemannian structure such that $\alpha$ is inner product preserving.
\ethm
{\it Proof}:\\
Fix a Fr\'echet-dense unital $\ast$-subalgebra $\cla$ over which $\alpha$ is algebraic. 
For proving the  implication $(ii)\raro(i)$, we consider $F=\alpha(f) d\alpha(g)-d\alpha(g)\alpha(f) \in \Omega^1(C^\infty(M)) \bar{\ot} \clq$
 and using (\ref{inner_prod_pres_111}), verify that $\lgl\lgl F, F\rgl\rgl=\alpha(\lgl\lgl(f dg -dg f), (fdg-dg f)\rgl\rgl)=0.$ This proves $F=0$, hence (i). 
 
The proof of the converse is basically an adaptation of arguments in  \cite{average_1}. We indicate only the steps where a modification of 
 the arguments of \cite{average_1} is necessary. We note from \cite{Huichi_huang} that 
 the antipode $\kappa$ (say) is defined and norm bounded on the whole of $\clq$.  We first observe  that the map $\Psi$ introduced  in Lemma 3.2 of 
 \cite{average_1} can be actually defined on a larger set, namely, for $F \in C^\infty(M) \ot \clq_0$. 
 Indeed, as the multiplication map $m$ is defined on $\clq \ot \clq_0$, we set $\Psi : C^\infty(M) \ot \clq_0 \raro C^\infty(M)$ 
  by 
	\begin{displaymath}
	\Psi(F)(x)=h \circ m \circ (\kappa \ot {\rm id})\left( (\alpha \ot {\rm id})(F)(x)\right),
	\end{displaymath}
	where $h$ denotes the Haar state. 
  Then the proof of the complete positivity of $\Psi$ as in Lemma 3.2 of \cite{average_1} goes through verbatim. Next, 
   choose  any Riemannian metric on $M$ and denote the corresponding $C^\infty(M)$-valued inner product 
   by $\lgl\lgl \cdot, \cdot \rgl\rgl$.
   Define \begin{displaymath} k(f,g):=\Psi(\lgl\lgl df_{(0)},dg_{(0)}\rgl\rgl \otimes f_{(1)}^{\ast}g_{(1)}),\end{displaymath} for $f,g \in \cla.$ Arguing along the lines of the proof of Lemma 
  3.3 of \cite{average_1} we can prove $\Psi(F \alpha(f))=\Psi(F)f$ for $F \in C^\infty(M) \ot \clq_0, f \in \cla$. From this, it follows that  
   $k(f,gh)=k(f,g)h+gk(f,h)$, i.e. for any fixed $f$, the map $g \mapsto k(f,g)$ is a derivation of $\cla$. As $\cla$ is Fr\'echet dense in $C^\infty(M)$, for any $x \in M$ 
    we can find $f_1, \ldots, f_n$ from $\cla$ such that $(f_1, \ldots, f_n)$ is a set  of  local coordinates in a neighbourhood $U$ (say) of $x$. The derivation 
     $k_f(\cdot)$ of the algebra generated by the coordinates  $f_1, \ldots, f_n$ must be of the form $k(f,\cdot)=Y^{U, f_1, \ldots, f_n}_f(\cdot)$, where $Y^{U, f_1, \ldots, f_n}_f$ 
     a vector field  defined and smooth on $U$. The local vector fields $Y^{U, f_1, \ldots, f_n}_f$'s do patch up consistently to give a smooth globally defined vector field 
     $Y_f$ on $M$. Indeed, 
      given two such local coordinates $(U, (f_1, \ldots, f_n))$ and $(V, (g_1, \ldots, g_n))$ with $U \bigcap V $ nonempty and $f_i, g_j \in \cla~ \forall i,j$, 
             we have \begin{displaymath}Y_f^{U, f_1, \ldots, f_n}(\phi)(x)=Y_f^{V, g_1, \ldots, g_n}(\phi)(x)=k(f,\phi)(x)\end{displaymath} for all $\phi \in \cla, x\in U \bigcap V$. 
             By the Fr\'echet density of $\cla$ and the obvious Fr\'echet continuity 
        of the locally smooth vector fields, they agree for all $\phi \in C^\infty(M)$, i.e. $Y_f^{U,f_1, \ldots, f_n}=Y_f^{V, g_1, \ldots, g_n}$ on $U\bigcap V$. 
        
        Clearly, for any $\eta \in \Omega^1(\cla)$, the map $\cla \ni f \mapsto \eta(Y_{\overline{f}})$ is a derivation, and arguing as before, we get a 
        globally defined  smooth vector field  $Z_\eta$ (say) such that $\eta(Y_{\overline{f}})=Z_\eta(f)=df(Z_\eta)$. 
         Define a sesquilinear form on $\Omega^1(M)$ by $\lgl\lgl \omega, \eta\rgl\rgl^\prime:=\overline{\omega}(Z_\eta)$. Clearly,
          \begin{displaymath}\lgl\lgl df, dg \rgl\rgl^\prime=d\overline{f}(Z_{dg})=dg(Y_f)=Y_f(g)=k(f,g)\end{displaymath} for $f,g \in \cla$. The rest of the arguments for proving that
           $\{ \lgl\lgl \cdot, \cdot \rgl\rgl_x^\prime$, $x \in M\}$ indeed gives the required invariant 
           Riemannian metric, are very similar to  those in \cite{average_1}. 
       In particular, the proof of Lemma 3.4  of \cite{average_1} goes through verbatim and  in the proof of Theorem 3.1 of \cite{average_1}, it is easy to observe  
       that  the assumption $\lgl\lgl df_{(0)},dg_{(0)}\rgl\rgl\ot f_{(1)}^{\ast}g_{(1)}
        \in \cla \ot \clq_0$ is not really necessary, as  the arguments go through even if $\lgl\lgl df_{(0)},dg_{(0)}\rgl\rgl\ot f_{(1)}^*g_{(1)}
        \in C^\infty(M) \ot \clq_0$.
        \qed
   \subsection{Lift of inner product preserving actions to bundles of one and higher forms}
From now on throughout Sections 3 and 4, we fix a smooth, faithful action $\alpha$ on $M$ and a Fr\'echet-dense unital $\ast$-subalgebra $\cla$ on which $\alpha$ is algebraic. 
In the present subsection we also assume that $\alpha$ preserves some Riemannian inner product.\\ 
\indent Using  the condition (\ref{inner_prod_pres_111}), we easily see that, whenever $\omega=\sum_{i=1}^l f_i dg_i=0$, where $f_i,g_i \in \cla$, we have 
$\sum_i \alpha(f_i)(d \ot {\rm id})(\alpha(g_i))
=0$. Thus, we get a well-defined map $\omega \mapsto d\alpha_{(1)}(\omega):= \sum_i \alpha(f_i)(d \ot {\rm id})(\alpha(g_i))$ from $\Omega^1(\cla)$ to $\Omega^1(\cla) \otimes \clq_0$.
Moreover, this has the properties  that $$ d\alpha_{(1)}(\omega f)=d\alpha_{(1)}(\omega) \alpha(f)=\alpha(f) d\alpha_{(1)}(\omega),~~ 
\lgl\lgl d\alpha_{(1)}(\omega), d\alpha_{(1)}(\eta)\rgl\rgl=\alpha(\lgl\omega, \eta\rgl)$$ for all $\omega, \eta \in \Omega^1(\cla), f \in \cla$. By Fr\'echet continuity of $\alpha$, 
we can 
 extend $d\alpha_{(1)}$ to  a continuous $\IC$-linear map from $\Omega^1(C^\infty(M))$ to $\Omega^1(C^\infty(M)) \bar{\ot} \clq$ satisfying similar properties. This motivates the 
  following definition.
  
\bdfn
\label{equiv_def}
Let $\cle$ be the module of smooth sections of a smooth hermitian bundle $E$  over $M$. 
A $\mathbb{C}$-linear map $\Gamma:\cle\raro \cle\bar{\ot}\clq$
is
said to be an $\alpha$ equivariant unitary representation of $\clq$ on 
$\cle$ if for all $\xi, \xi^\prime \in \cle,~ f\in C^\infty(M),$\\ 
1. $\Gamma(\xi f)=\Gamma(\xi)\alpha(f)=\alpha(f) \Gamma(\xi)$,\\
2. $\lgl\lgl\Gamma(\xi),\Gamma(\xi^{\prime})\rgl\rgl=\alpha(\lgl\lgl\xi,\xi^{\prime}\rgl\rgl)$,\\
3. $((\Gamma \overline{\ot} {\rm id})\Gamma) (\xi)(m)= \Delta(\Gamma(\xi)(m))$ $\forall \xi \in \cle,~m \in M$ (co-associativity),\\
4. The right $\clq$-linear  span of $ \Gamma(\cle)$ is Fr\'echet-dense in $\cle\bar{\ot}\clq$ (non
degeneracy).
\edfn
Note that the condition 2. in the definition above  allows one to define
$(\Gamma\overline{\ot} {\rm id})$ as the extension of $\Gamma \overline{\ot} {\rm id}$ to the completion of $\cle \ot \clq.$ We have also used the identification 
 of $\Gamma(\xi)$ as a $\clq$-valued smooth section. Moreover, 3. can be interpreted (formally) as
 $(\Gamma \overline{\ot} {\rm id})\Gamma=({\rm id} \ot \Delta)$. 
  We also simply say `equivariant unitary representation' if the action $\alpha$ is understood from the context. 
 In particular,  when $\cle$ is the trivial $C^\infty(M)$-bimodule of rank $N$,
 we have the following:
\blmma
\label{unitary}
Given an $\alpha$ equivariant unitary representation $\Gamma$ of $\clq$ on
$\mathbb{C}^{N}\ot C^\infty(M)$ such that $\Gamma(e_{i}\ot
1)=\sum_{j=1}^{N}e_{j}\ot b_{ji}$, $b_{ij}\in C^\infty(M, \clq)$ for
all $i,j=1,...,N$, where $\{e_{i};i=1,...,N \}$ is an orthonormal basis of
$\mathbb{C}^{N}$, then $B=((b_{ij}))_{i,j=1,....,N}$ is a unitary element of
$M_{N}(C^\infty(M, \clq))$.
\elmma
{\it Proof}:\\
Clearly, $B$ can also be thought of  as an element of the $C^*$ algebra $M_N(\clc),$ where $\clc:=C(M) \hat{\ot} \clq$, 
which is isomorphic with the set of right $\clc$-linear, adjointable maps 
 on the Hilbert module $\clf:={\mathbb C}^N \ot \clc$.  
Viewing $B$ as a (right $\clc$-linear) map on $\clf$, we first claim that $B^*B=I_\clf$. Indeed, 
 $ \lgl\lgl\Gamma(e_{i}\ot 1),\Gamma(e_{j}\ot 1)\rgl\rgl=\alpha(\lgl\lgl e_{i}\ot 1,e_{j}\ot
1\rgl\rgl)$, which implies 
 \begin{displaymath}\sum_{k=1}^{N}b_{ki}^{\ast}b_{kj}= \alpha(\delta_{ij}1)=\delta_{ij}1_{\clc},\end{displaymath}  proving   the claim. 
It now suffices to prove that the  range of $B$ is
dense in $\clf$. To this end, observe that
 \begin{eqnarray*}
 \Gamma(e_{i}\ot f)q &=& \sum_{j}e_{j}\ot b_{ji}\alpha(f)(1\ot q)\\
&=& B(e_{i}\ot 1_{\clc})\alpha(f)(1\ot q)\subseteq {\rm Im}(B),
\end{eqnarray*}
 for $f\in C^\infty(M),q\in\clq$. Thus, ${\rm Im}(B)$ contains the right $\clq$-linear span of $\Gamma({\mathbb C}^N \ot C^\infty(M))$, which is 
Fre\'chet-dense  in ${\mathbb C}^N \ot C^\infty(M, \clq)$, hence norm-dense in $\clf$ as well. \qed
 
 It is clear that  $d \alpha_{(1)}$ is an equivariant unitary representation on $\Omega^1(C^\infty(M))$. 
 We now go further to construct similar representation on higher forms. 
\blmma
\label{lift_eq_rep}
For each $k \geq 1$, there is an equivariant unitary representation $U^{(k)}$ of $\clq$ on each $\Omega^k(C^\infty(M))$, 
  given on the dense subspace $\Omega^k(\cla)$  by $$U^{(k)}(\omega_1 \ot_\cla \omega_2 \ot_\cla \ldots \ot_\cla \omega_k)=\omega_{1(0)} 
 \ot_\cla...\ot_{\cla} \omega_{k(0)}\ot  \omega_{1(1)}\omega_{2(1)} \ldots \omega_{k(1)},$$
 where $\omega_i \in \Omega^1(\cla)~\forall i$ and the Sweedler type notation for comodule maps  has been used. 
 \elmma
 {\it Proof:}\\ For $k=1$, $U^{(1)}=d\alpha_{(1)}$ and there is nothing to prove. Let us prove the result for $k=2$ only, as the proof for $k \geq 3$ is  very similar. 
 Using the condition 2. of Definition \ref{equiv_def} for $d\alpha_{(1)}$, we verify that $U^{(2)}$ is well defined.  
  Indeed, for $\omega, \eta \in \Omega^1(\cla)$, $f \in \cla$, we have \begin{eqnarray*}(\omega f)_{(0)} \ot_\cla \eta_{(0)} \ot (\omega f)_{(1)}\eta_{(1)}
  &=&\omega_{(0)} f_{(0)}\ot_\cla \eta_{(0)} \ot \omega_{(1)} f_{(1)} \eta_{(1)}\\
  &=& \omega_{(0)} \ot_\cla f_{(0)} \eta_{(0)} \ot \omega_{(1)} f_{(1)} \eta_{(1)}\\
  &=&\omega_{(0)}\ot_\cla (f \eta)_{(0)}  \ot \omega_{(1)}(f \eta)_{(1)}.\end{eqnarray*} Next, $U^{(2)} $ on $\Omega^2(\cla)$ is clearly an algebraic (co)-representation of $\clq_0$, 
   so it satisfies the co-associativity condition 3. and also $({\rm id} \ot \epsilon) \circ U^{(2)}={\rm id}$. Moreover, using conditions 1. and 2. of Definition \ref{equiv_def} 
   for $d\alpha_{(1)}$ and for 
   $\omega,\eta, \omega^\prime,\eta^\prime \in \Omega^1(\cla), $ we have for $x\in M$, 
  \begin{eqnarray*} 
	\lgl\lgl U^{(2)}(\omega \ot_\cla \eta), U^{(2)}(\omega^\prime \ot_\cla \eta^\prime)\rgl\rgl(x)
   &=&\lgl\lgl\omega_{(0)},\omega^{\prime}_{(0)}\rgl\rgl(x)\lgl\lgl\eta_{(0)}, \eta^{\prime}_{(0)}\rgl\rgl(x)\eta_{(1)}^{\ast}\omega_{(1)}^{\ast}\omega_{(1)}^{\prime}\eta^{\prime}_{(1)}\\
   &=&\lgl\lgl\eta_{(0)},\eta^{\prime}_{(0)}\rgl\rgl(x)\eta^{\ast}_{(1)}(\lgl\lgl\omega_{(0)},\omega^{\prime}_{(0)}\rgl\rgl(x)\omega^{\ast}_{(1)}\omega^{\prime}_{(1)})\eta_{(1)}^{\prime}\\
   &=&\lgl\lgl\eta_{(0)},\eta^{\prime}_{(0)}\rgl\rgl(x)\eta^{\ast}_{(1)}\alpha(\lgl\lgl\omega,\omega^{\prime}\rgl\rgl)(x)\eta^{\prime}_{(1)}
	\end{eqnarray*}
	which,  
   by (\ref{inner_prod_pres_111}), becomes 
   \begin{eqnarray*}
	\alpha(\lgl\lgl\omega,\omega^{\prime}\rgl\rgl)(x)\lgl\lgl\eta_{(0)},\eta^{\prime}_{(0)}\rgl\rgl(x)\eta^{\ast}_{(1)}\eta_{(1)}^{\prime}
   &=&\alpha(\lgl\lgl\omega,\omega^{\prime}\rgl\rgl)(x)\alpha(\lgl\lgl\eta,\eta^{\prime}\rgl\rgl)(x)\\
   &=&\alpha(\lgl\lgl \omega \ot_\cla \eta, \omega^\prime \ot_\cla \eta^\prime\rgl\rgl)(x).
	\end{eqnarray*}
   This proves condition 2. for $U^{(2)}$ on $\Omega^2(\cla)$,
     hence $U^{(2)}$ extends to $\Omega^2(C^\infty(M))$ by continuity, and the extension is easily seen to satisfy the conditions 1, 2 and 3. 
     Finally, as $U^{(2)}(\Omega^2(\cla))\clq_0$ is total in $\Omega^2(\cla) \ot \clq_0,$ by continuity we get the density condition of 4. 
     This completes the proof. \qed
    
    Recall the subalgebra $\clq_x$ defined by (\ref{def_q_x}) in Subsection 3.1.
 \blmma
  \label{commute}
  For $x \in M$, $\clq_x$ is commutative. 
   \elmma
   {\it Proof:}\\
   Recall the decomposition of $\Omega^2(\cla)=\clf_s(\cla) \oplus \Lambda^{2}(\cla)$ and the description of $\clf_s(\cla)$ given by Proposition \ref{symm_description}. 
  Let  $f_i,g_i \in \cla$ be such that $\sum_i f_idg_i=0$. By applying $d\alpha_{(1)}$ we get 
    $\sum _i f_{i(0)}dg_{i(0)}\ot f_{i(1)}g_{i(1)}=0$, hence $\sum f_{i(0)}dg_{(i)0}\phi(f_{i_{(1)}}g_{i_{(1)}})=0$  for every linear functional $\phi$ on $\clq_0$. So we have 
    \begin{displaymath}({\rm id} \ot \phi)(U^{(2)}(\sum_i df_i \ot_\cla dg_i))
   =\sum_i df_{i(0)} \ot_\cla dg_{i(0)} \phi( f_{i(1)}g_{i(1)}) \in \clf_s(\cla).
	\end{displaymath}
   Thus, $U^{(2)}(\clf_s(\cla)) \subseteq \clf_s(\cla) \ot \clq_0$. Moreover, if $\omega \in \Lambda^2(\cla)=\clf_s(\cla)^\perp$, we have \begin{displaymath}\lgl\lgl U(\omega), U(\eta)q\rgl\rgl=
   \alpha(\lgl\lgl \omega, \eta\rgl\rgl)q=0\end{displaymath} for any $\eta \in \clf_s(\cla), q\in \clq_0$. But  $U^{(2)}|_{\clf_s(\cla)}$ is an algebraic (co)-representation which implies that
    $U^{(2)}(\clf_s(\cla))\clq_0=\clf_s \ot \clq_0$. In particular, $\eta \ot 1 \in U^{(2)}(\clf_s(\cla))\clq_0$, which implies
     \begin{displaymath}\lgl\lgl\eta,({\rm id} \ot \phi)(U^{(2)}(\omega))\rgl\rgl=({\rm id}\ot\phi)(\lgl\lgl(\eta \ot 1),U^{(2)}(\omega)\rgl\rgl)=0\end{displaymath} for $\omega \in \Lambda^2(\cla), \eta \in \clf_s(\cla)$, $\phi \in \clq_0^\prime$.
      Hence $U^{(2)}(\Lambda^2(\cla)) \subseteq \Lambda^2(\cla) \ot \clq_0$ and the restriction of $U^{(2)}$ to 
   $\Lambda^2(\cla)$ gives a unitary equivariant representation of $\clq_{0}$ on $\Lambda^2(\cla)$.
    

Now, choose smooth one-forms $\{
\omega_1,\ldots, \omega_n  \}$ such that they form a basis of $T^*M$ at $x$. Recall the notation $X_{\omega}$ from Subsection \ref{form}. Clearly, it is enough to 
 prove that  $F_i(x):=X_{\omega_i}(\alpha(f))(x)$ and $G_j(x):=X_{\omega_j}(\alpha(g))(x)$ commute for $f,g \in C^\infty(M)$ and $\forall i,j=1, \ldots, n$.  
 Indeed,  $d \alpha(f)(x)=\sum_i \omega_i(x) F_i(x), $ 
$d \alpha(g)(x)=\sum_i \omega_i(x) G_i(x)$. Now $U^{(2)} $ leaves
invariant the submodules of
symmetric and antisymmetric tensor product of $\Lambda^1(\cla)$, 
in
particular, $C^s_{ij}=C^s_{ji}, $ $C^a_{ij}=-C^a_{ji}$ for all $i,j$, where
$C^s_{ij}$ and $C^a_{ij}$  denote the $\clq$-valued coefficient of $w_i(x)\ot
w_j(x)$ in the expression of $U^{(2)}(df \ot dg+dg \ot df)|_{x}$ and
$U^{(2)}(df \ot dg - dg \ot df)|_{x}$ respectively. By a simple
calculation using these relations, we get the commutativity of $F_i(x), G_j(x)$
for all $i,j$.
\qed\vspace{2mm}\\    
 We denote the following $\ast$-subalgebra of $C^\infty(M, \clq)$ by $\clc$, which is commutative by Lemma \ref{commute}:
 $$ \clc :=\{ F :~F(x) \in \clq_x~\forall x \}.$$ 
  Viewing $d\alpha_{(1)}(\omega)$ as a $\clq$-valued smooth section, let us write $d\alpha_{(1)}(\omega)(m)=d\omega_{(0)}(m) \ot  \omega_{(1)} \in T^*_m M \ot \clq_0 $ 
  for $\omega \in \Omega^1(\cla)$, $m \in M$. Then it follows 
   from the above lemma that $\forall \omega, \eta \in \Omega^1(\cla)$, \be \label{comm111}
   \omega_{(0)}(m)\ot \eta_{(0)}(m)\ot  \omega_{(1)}\eta_{(1)}=\omega_{(0)}(m) \ot \eta_{(0)}(m) \ot \eta_{(1)}\omega_{(1)}.
   \ee
   
\bppsn
\label{lifting_rep}
If $\alpha$ is a Riemannian inner product preserving action, then for every $k\geq 0$, there exist an $\alpha$-equivariant representation $d\alpha_{(k)}$ of $\clq$ on $\Lambda^{k}(C^{\infty}(M))$ given by
$$d\alpha_{(k)}(f_{0}df_{1}\wedge...\wedge df_{k})=\alpha(f_{0})d(\alpha(f_{1}))\wedge...\wedge d(\alpha(f_{k})).$$ Moreover, we have 
\be \label{777} d\alpha_{(k+l)}(\omega \wedge \eta)=d\alpha_{(k)}(\omega) \wedge d\alpha_{(l)}(\eta)\ee for $\omega \in \Lambda^k(M, \clq)$, $\eta \in \Lambda^l(M, \clq)$. Here we have used 
 the wedge product between $\clq$-valued forms introduced earlier. 
\eppsn
{\it Proof:}\\ It follows from (\ref{comm111})  that $U^{(k)}$ commutes with the action of $S_k$ on $\Omega^k(C^\infty(M))$, hence leaves each of 
 the spectral subspaces for the $S_k$-action invariant, in particular the range of $P_{\rm sgn}$, i.e. $\Lambda^k(C^\infty(M))$. We take the restriction of $U^{(k)}$ 
  on $\Lambda^k(C^\infty(M))$ as the definition of $d\alpha_{(k)}$.  
  
  It is clear from the definition that $d\alpha_{(k)}(f\omega)=\alpha(f)d\alpha_{(k)}(\omega)$. We can now 
   prove (\ref{777}) by considering $\omega=fdf_1\wedge \ldots \wedge df_k, $ $\eta=gdg_1 \wedge \ldots \wedge dg_l$ and using the commutativity of $\clq_m$ $\forall m$.\qed

  \subsection{Lift to the orthonormal frame bundle}
  Let $O(M)$ be the bundle of orthonormal frames of $M$. 
    We can identify this as  a subbundle of the direct
sum of $n$ copies of cotangent bundle, say
$E=\underbrace{T^{\ast}(M)_{\mathbb{R}}\oplus...\oplus T^{\ast}(M)_{\mathbb{R}}}_{n-copies}$. Let $\overline{\xi}$ denote the complex conjugate of a complexified  
co-tangent vector $\xi \in T^*_m M$, and for a (complexified) one-form $\omega$, define $\overline{\omega}(m)=\overline{\omega(m)}$ for all $m \in M$. 
Indeed, $$ O(M)=\{(m,\underline{\omega}):~m \in M,~\underline{\omega}=(\omega_1,\ldots, \omega_n),\omega_i \in T^{\ast}_{m}M,\overline{\omega_{i}}=\omega_{i}\forall i~<\omega_{i},\omega_{j}>=\delta_{ij}\}.$$ 
 Denote by $\pi$ the projection for both the bundles $E$ and $O(M)$.  The fibre of $O(M)$ at
each point is $K=\{\underline{v}\equiv(v_{1},...,v_{n}): v_{i}\in\mathbb{R}^{n}
\forall i \ {\rm and} \ <v_{i},v_{j}>=\delta_{ij} \}$. We identify $\underline{v} \in K$ with the  $n\times n$ real orthogonal matrix $((v_{ik}))$, where 
 $v_i=(v_{i1}, \ldots, v_{in})$.  Thus $K \cong O_n(\IR)$.  Let $\cle$ be the complexification of the space of smooth sections of $E$, which is 
 $\cle=\Omega^1(C^\infty(M)) \oplus \ldots \oplus \Omega^1(C^\infty(M))$
     ($n$ copies).  
 Clearly,  $\cle_0:=\Omega^{1}(\cla)\oplus \ldots \Omega^1(\cla)$ ($n$ copies)  is Fre\'chet-dense in $\cle$. 
For any subset $U $ of $M$ we denote by $E_U$ and $O(U)$ the restrictions of the bundles $E$ and $O(M)$ to $U$ respectively.

  There is a
natural $C^{\infty}(M)$ valued inner product $\lgl\lgl,\rgl\rgl^{\cle}$ on $\cle$,
given by $\lgl\lgl\underline{\omega},\underline{\eta}\rgl\rgl^{\cle}(m):=
\sum_{i=1}^{n}\lgl\omega_{i}(m),\eta_{i}(m)\rgl$. 
  For a local coordinate chart $(U,(x_{1},x_{2},...,x_{n})),$ a choice of real one forms $\{\omega_{1}^{\prime},\omega_{2}^{\prime},...,\omega_{n}^{\prime}\}$, i.e. $\overline{\omega_{i}}=\omega_{i}$ for all $i$, 
  will be  called {\it $U$-orthonormal } if for all $x\in U$, $\{\omega_{1}^{\prime}(x),\omega_{2}^{\prime}(x),...,\omega_{n}^{\prime}(x)\}$ are orthonormal with respect to 
  the Riemannian metric. Let us fix some $U$-orthonormal set $\{\omega^{\prime}_{1},...,\omega^{\prime}_{n}\}$. For $\xi \in \Omega^{1}(C^\infty(M))$, $1 \leq i \leq n$, let  $\theta^{i}_{\xi} \in
C^\infty(O(M))$
 be given by $ \theta^{i}_{\xi}(e):=\lgl \lgl{\omega_{i}}, \xi \rgl \rgl  
(\pi(e)),$ where $e=(m,\underline{\omega})$ . Define $s^{U}_{ij}\equiv s^{(U,\underline{\omega^{\prime}})}_{ij}\in C^{\infty}(E)$ by $s^{(U,\underline{\omega^{\prime}})}_{ij}(m,\underline{\omega})=<\omega_{i}(m),\omega^{\prime}_{j}(m)>$. Let $t^U_{ij}\equiv t^{U,\underline{\omega^\prime}}_{ij}$ be the restriction of $s^U_{ij}$ to $O(M)$. Then $t^{U}_{ij}=\theta^{i}_{\omega^{\prime}_{j}}$.   
Note that $d\alpha_{(n)}(\overline{\omega})=\overline{d \alpha_{(n)}(\omega)}$ $\forall \omega \in \Omega^n(C^\infty(M))$, which can be verified by writing $\omega= f_{0}df_{1}\wedge df_{2}\wedge...\wedge df_{n}$ and using $\alpha(\overline{f_{i}})=(\alpha(f_{i}))^{\ast}$ for all $i=0,...,n$. Hence $s^{U}_{ij}$ and $t_{ij}^{U}$ are self adjoint.\\ 
\indent From the 
 definition of $O(M)$, we observe that $(( t^U_{ij}(e) )) \in O_n(\IR)$ $\forall e \in O(M)$. 
Clearly, $\{ s^U_{ij}, i,j=1,\ldots, n \}$ forms a set of local coordinate functions in the fibre direction for $\pi^{-1}(U)$, hence 
 the algebra generated by the functions of the form $(f \circ \pi) s^U_{ij}$, with $f$ smooth and having support in a compact subset of $U$, constitute a Fr\'echet-dense subset of 
 $C_c^\infty(E_U)$. By restriction, we get a similar dense algebra, say $\clp^\infty_U$, of $C_c^\infty(O(U))$,
  generated by $(f \circ \pi) t^U_{ij}$'s, $f \in C_c^\infty(U)$.  Moreover, $(f \circ \pi)t^U_{ij}=\theta^i_{f\omega^{\prime}_j}$, which shows that the algebra generated 
  by $\theta^i_\omega$,
   with $i=1, \ldots, n$ and $\omega \in \Omega^1(C^\infty(M))$, is Fr\'echet dense in $C^\infty(O(M))$ too.

Define $T^{U,\underline{\omega^\prime}}_{ij} \equiv T^{U}_{ij} \in C^\infty(O(M), \clq)$ by:
$$T^{U}_{ij}(m,\underline{\omega}):=\lgl\lgl(\omega_{i} \ot
1_\clq),d\alpha_{(1)}(\omega^{\prime}_{j})\rgl\rgl^{\cle\bar{\ot}\clq}(m),$$ 
   where $\lgl\lgl,\rgl\rgl^{\cle\bar{\ot}\clq}$ denotes the $C^{\infty}(M,\clq)$ valued
inner product as before
and $\underline{\omega}=(\omega_{1},...,\omega_{n})$.
   Clearly,  $T^{U}_{ij}(e)\in \clq_{\pi(e)}~\forall e \in O(M), \forall i,j=1, \ldots n,$ hence $T_{ij}$'s commute among themselves.

Let us fix a coordinate neighbourhood $U$, $U$-orthonormal one forms $\omega^\prime_j$'s and consider the corresponding $t^{U}_{ij}, T^{U}_{ij}$'s.
 
     \blmma 
     \label{orthogonality}
     For any smooth real-valued function $\chi$ supported in $U$, we have
     \be \label{orthogonality_T} (\alpha(\chi)\circ \pi) \sum_j T^U_{ij}T^U_{lj}=
     (\alpha(\chi)\circ \pi)\delta_{il},\ee for all
$i,l=1,\ldots,n$.
 
    \elmma
    {\it Proof:}\\ Fix $e =(m, (\omega_1,\ldots, \omega_n)) \in O(M)$, $m=\pi(e)$. Let $\gamma$ be a character (multiplicative linear functional) on the commutative 
     $C^*$ algebra $\clq_m$ and $u_j:=({\rm id} \ot \gamma)(d\alpha_{(1)}(\omega^\prime_j)(m)) \in T^{\ast}_{m} M$. 
      By a simple calculation using  the facts that $d\alpha_{(1)}$ is inner-product preserving and ${\omega^\prime}_j$'s form an orthonormal basis of $T^{*}M$ at every point in the 
       support of $\chi$,  we obtain 
       \bean &&\lefteqn{\gamma(\alpha(\chi)(m))^2 \lgl u_i,u_l\rgl}\\
       &=& \gamma \left( \lgl\lgl d\alpha_{(1)}(\chi {\omega^\prime}_i), d\alpha_{(1)} (\chi
{\omega^\prime}_l)\rgl\rgl^{\cle\bar{\ot}\clq}(m) \right)\\
       &=& \gamma \left( \alpha(\lgl\lgl\chi{\omega^\prime}_i,
\chi{\omega^\prime}_l\rgl\rgl^{\cle})(m)\right)\\
       &=& \gamma(\alpha(\chi^2)(m))\delta_{il},\\
       \eean
     where $\lgl,\rgl$ denotes the inner product of $T_{m}^{\ast}(M)$.
       Thus, in case $\gamma(\alpha(\chi)(m)))$ is nonzero, $u_1,\ldots, u_n$ is an orthonormal basis of $T^{*}_{m}M$ and by multiplicativity of 
       $\gamma$ and self-adjointness of $T^U_{ij}$'s it is easy to see that
$\gamma(T^U_{ij}(e) )=\lgl\omega_i(\pi(e)),u_j\rgl=\lgl u_j,\omega_i(\pi(e))\rgl$. We have
                    \bean &&\lefteqn{ \gamma \left( \alpha(\chi^2)(\pi(e))  \sum_j T^U_{ij}(e)T^U_{lj}(e) \right)}\\
   &=& \gamma(\alpha(\chi^2)(\pi(e))) \sum_j \lgl\omega_i(\pi(e)), u_j\rgl\lgl u_j,\omega_l(\pi(e))\rgl\\
   &=& \gamma(\alpha(\chi^2)(\pi(e))) \lgl\omega_i(\pi(e)), \omega_l(\pi(e))\rgl= \delta_{il}\gamma(\alpha(\chi^2)(\pi(e))).\\
   \eean
   This equality is also true trivially when $\gamma(\alpha(\chi)(\pi(e)))=0$. Hence  
    $(\alpha(\chi)  (\pi(e)))^{2}X(e)= \delta_{il} (\alpha(\chi)(\pi(e)))^{2}$, where
    $X=\sum_{j}T^U_{ij}T^U_{lj}$.  As
    $\alpha(\chi)(\pi(e))$ is self adjoint, the equality (\ref{orthogonality_T} ) follows.
            \qed

    \blmma 
    \label{basis_change}
    Let $U,V$ be two coordinate neighbourhoods, $U \bigcap V\neq\emptyset$. Also let $\{ \omega^\prime_1, \ldots, \omega^\prime_n\}$ and 
    $\{ \omega^{\prime \prime}_1, \ldots, \omega^{\prime \prime}_n\}$ be $U$-orthonormal and $V$-orthonormal one-forms respectively. Then there are smooth functions $f_{jk}$ such that 
    for every $f \in C_c^\infty(U \bigcap V)$, we have 
    \be \label{basis_1}(f \circ \pi) t^{V, \underline{\omega^{\prime \prime}}}_{ij}=
                  \sum_k (f \circ \pi) (f_{jk} \circ \pi) t^{V, \underline{\omega^\prime}}_{ik},\ee
                  \be \label{basis_2}
                  (\alpha(f)\circ \pi)T^{V, \underline{\omega^{\prime \prime}}}_{ij}=\sum_k (\alpha(ff_{jk}) \circ \pi)T^{V, \underline{\omega^{\prime}}}_{ik}.\ee
    \elmma
 {\it Proof :}\\ 
  As both $\omega^\prime_i$'s and $\omega^{\prime \prime}_i$'s are bases of $T^{*}_{m} M$ at each $m \in U$, there are smooth functions $f_{jk}$ such that
                  \be \label{compare} f \omega^{\prime \prime}_j=\sum_k f  f_{jk} \omega^\prime_k\ee for any $f \in C_c^\infty(U)$. This implies (\ref{basis_1}). 
                 We obtain (\ref{basis_2}) by applying  $d\alpha_{(1)}$ on both sides of (\ref{compare}).\qed

   \blmma
   \label{density} 
   There is a well-defined, norm-contractive, $\ast$-homomorphism $\eta_U : C_{c}^{\infty}(U) \raro C^\infty(O(M), \clq)$ satisfying 
   \be \label{formula_1} \eta_U((f\circ \pi)t^U_{ij})=(\alpha(f)\circ \pi)
   T^U_{ij}.\ee  for all $f \in C_c^\infty(U).$ Moreover, \\
   (i) $\eta_U$ is continuous w.r.t. the Fr\'echet topology.\\
   (ii) $\eta_U$ does not depend on the choice of 
   the sections $\omega^\prime_i$'s.\\
   (iii) For two coordinate neighbourhoods $U$ and $V$, with $U \bigcap V$  nonempty and for $F \in C_c^\infty(O(U))$, $G \in C_c^\infty(O(V)),$
    we have $\eta_U(FG)=\eta_V(FG)=\eta_U(F) \eta_V(G)$.
      \elmma
      
     {\it Proof:}\\
     The map $\phi : O(U) \raro U \times K \cong U \times O_n(\IR)$ given by $\phi(e)=(\pi(e), ((t^U_{ij}))_{i,j=1, \ldots n})$ 
     is a diffeomorphism. For $F \in C_c^\infty(O(U))$,
     and a given $\underline{v} =(v_1, \ldots, v_n) \in K$, let $F_{\underline{v}} \in C_c^\infty(U) \subset C^\infty(M)$ 
     given by $F_{\underline{v}}(m)=(F \circ \phi^{-1})(m, \underline{v}),$ 
     $m \in U$ and $0$ for $m \not\in U$.
          Let $H_m(F) \in C^\infty(K, \clq_m) \subseteq C^\infty(K, \clq) $ be defined by $H_m(F)(\underline{v})=\alpha(F_{\underline{v}})(m), m \in M$. Clearly, 
          as the maps $\underline{v}
          \mapsto F_{\underline{v}}$ and $\alpha$ are Fr\'echet continuous, $\underline{v} \mapsto H_m(F_{\underline{v}})$ is Fr\'echet continuous too 
           for any fixed $m$. We want to define an element $\eta_U(F)(e) \in \cle_{\pi(e)}$ by specifying $\eta_U(F)(e)(\gamma) \equiv \gamma(\eta_U(F)(e))$ 
           for $\gamma \in \hat{\clq}_{\pi(e)}$,
           where $\hat{\clq}_{\pi(e)}$ denotes the set of  characters 
            the commutative $C^{\ast}$ algebra $\clq_{\pi(e)} \cong C(\hat{\clq}_{\pi(e)})$, with the weak $\ast$ topology.
              Let $\chi \in C_c^\infty(U)$ be such that $0 \leq \chi \leq 1$ and $\chi \equiv 1$ on $\pi({\rm Supp}(F))$,  so that $(\chi \circ \pi) F=F$. Define 
             $\gamma(\eta_U(F)(e))=0$ if $\gamma(\alpha(\chi)(\pi(e)))=0$. Otherwise, by Lemma \ref{orthogonality}, $\underline{\tau(e)}:=(\gamma(T^U_{ij}(e)),i,j=1, \ldots n)$ 
             belongs to $K$, hence it makes sense
              to define the following: 
							\begin{displaymath}
							\gamma(\eta_U(F)(e)):=\gamma(\alpha(\chi)(\pi(e))) \gamma (H_{\pi(e)}(F)(\underline{\tau(e)})). 
               \end{displaymath}                      
              It follows from the continuity of $H_m(F)$ that $\gamma \mapsto \gamma(\eta_U(F)(e))$ is continuous and as
              $\alpha$ and $\gamma$ are $\ast$-homomorphisms, hence norm-contractive with
								\begin{displaymath}\| \gamma(\eta_U(F)(e))\| \leq \| \chi \|_\infty 
              \sup_{m \in M,~\underline{v} \in K} |F_{\underline{v}}(m)|_\infty \leq \| F\|_\infty.
							\end{displaymath}
							This allows us 
               to extend $\eta_U$ as a norm-contractive map from $C_c(E_U)$ to $C(E, \clq)$. 
               It is also easy to see  that $F \mapsto \gamma(\eta_U(F)(e))$ is $\ast$-homomorphic
               for every $\gamma$, which implies $C_c(E_U) \ni F \mapsto \eta_U(F) \in C(E, \clq)$ is indeed $\ast$-homomorphic. 
               To verify (\ref{formula_1}), it is enough to observe that, by definition, $H_m((f \circ \pi) t^U_{ij})(\underline{v})=\alpha(f)(m) v_{ij}$. 
               
               At this point, we claim that the  definition of $\eta_U$ does not really depend on the choice of $\chi$. Indeed, if $\chi^\prime$ is another smooth function supported in $U$,
               with $(\chi^\prime \circ \pi) F=F$, denoting by $\eta^\prime_U$ the analogue of $\eta_U$ using $\chi^\prime$ instead of $\chi$, 
                we'll have $F_{\underline{v}}=\chi F_{\underline{v}}=\chi^\prime F_{\underline{v}}$, 
                hence \be \label{eqn_1000} H_m(F)=\alpha(\chi)(m)H_m(F)=\alpha(\chi^\prime)(m)H_m(F).\ee For a given $e$ and $\gamma$, if $\gamma(\alpha(\chi)(\pi(e)))=0$ but $\gamma(\alpha(\chi^\prime)(\pi(e)))$
                 is nonzero, we have from (\ref{eqn_1000}) that \begin{displaymath}\gamma(H_{\pi(e)}(F)(\cdot))
                =0=\gamma(\alpha(\chi^\prime)(\pi(e)))\gamma(H_{\pi(e)}(F)(\cdot)),\end{displaymath} so that $\gamma(\eta_U(F)(e))=\gamma(\eta^\prime_U(F)(e))=0$. Otherwise also, 
                the equality of $\gamma(\eta_U(e))$ and $\gamma(\eta^\prime_U(e))$ follows from (\ref{eqn_1000}). This proves our claim.

                Next we verify that the definition is also independent of the choice of the one-forms $\omega^\prime_j, j=1, \ldots, n.$ Let $\omega^{\prime \prime}_j, j=1, \ldots, n$
                 be another such choice and let $\gamma^U_{ij}=t^{\underline{\omega^{\prime \prime}}, U}_{ij}$, $\Gamma^U_{ij}=T^{\underline{\omega^{\prime \prime}}, U}_{ij}$. 
                 Denote by $\zeta_U$ the analogue of $\eta_U$ defined using $\omega^{\prime \prime}_j$'s. Clearly, 
                  $\zeta_U$ is uniquely determined by $\zeta_U((f \circ \pi) \gamma^U_{ij})$ for all $f \in C_c^\infty(U)$, hence 
                   it suffices to prove $\eta_U((f \circ \pi) \gamma^U_{ij})=(\alpha(f) \circ \pi)\Gamma^U_{ij}.$ However, this follows from Lemma \ref{basis_change}, taking $U=V$, 
                    first by
                    applying $\eta_U$ on the expression of $\gamma^U_{ij}$ obtained from  (\ref{basis_1}) and then using (\ref{basis_2}).

                 We can prove (iii) also by applying Lemma \ref{basis_change}. In the notation of Lemma \ref{basis_change},  it is enough to show 
                 $$\eta_U((f\circ \pi)(g \circ \pi) t^{U, \underline{\omega^\prime}}_{ij} t^{V, \underline{\omega^{\prime \prime}}}_{kl})=
                 \eta_U((f\circ \pi)t^{U, \underline{\omega^\prime}}_{ij})\eta_V((g\circ \pi)t^{V, \underline{\omega^{\prime \prime}}}_{kl})$$
                 $\forall$ $f \in C_c^\infty(U)$, $g \in C_c^\infty(V).$ Substituting $t^{V, \underline{\omega^{\prime \prime}}}_{kl}$ by 
                 $t^{U, \underline{\omega^\prime}}_{pq}$'s using (\ref{basis_1}) of Lemma \ref{basis_change}, and the definition of $\eta_U$, the left 
               hand side is seen to be equal to  $(\alpha(f) \circ \pi)T^{U, \underline{\omega^\prime}}_{ij})\sum_p (\alpha(gf_{lp})\circ \pi )T^{U, \underline{\omega^\prime}}_{kp}$,
               which coincides with the right hand side by (\ref{basis_2}) of Lemma \ref{basis_change}, replacing $(ij)$ by $(kl)$.

               To prove the Fr\'echet continuity of the map at a point $e_0 \in U$, choose open set $U_0$ such that $e_0 \in U_0 \subset \overline{U_0} \subset U$. 
                Consider the embedding $\psi$ (say)  $e \mapsto (x_1, \ldots, x_n, t^U_{ij}, i,j=1, \ldots, n )$ of $U \times K  \subset \IR^{n+n^2}$.
                           As the map $(m,\underline{v}) \mapsto H(m,\underline{v}):= H_m(F)(\underline{v})$
               is  smooth 
                on  the compact set $\overline{U_0}  \times K $  we can choose an open 
                neighbourhood $W$ of $\psi(\overline{U_0} \times K) \subset \IR^{n+n^2}$ such that   $H \circ \psi^{-1}$ extends as a smooth $\clq$-valued function $H_1$ (say) on $W$.  
                Choose $\lambda \in C_c^\infty(\IR^{n+n^2})$ satisfying $\lambda(z)=1$
                 for all $z \in \psi(\overline{U_0} \times K)$, hence $\lambda H_1=H\circ \psi^{-1}$ on $\psi(\overline{U_0} \times K).$ 
                  For $e \in O(U)$, we claim  the following:
                \be  \label{eqn1}
                \eta(F)(e):=\alpha(\chi)(\pi(e)) \int
                \widehat{\lambda H_1}(\underline{\mu},\underline{\xi}){\rm exp}\left( -2 \pi i \sum_k  \mu_k x_k(\pi(e))  -2 \pi i \sum_{jl}\xi_{jl} T^U_{jl}(e) \right) d\underline{\mu} d\underline{\xi},\ee
                   where the integration is over $(\underline{\mu},\underline{\xi}) \in \IR^n\times\IR^{n^2}$ and $\widehat{\lambda H_1}$ denotes the Fourier transform of the smooth, compactly supported function $\lambda H_1$. It is clear from the above expression and the smoothness of $T_{jl}^{U}$ that $\eta(F)$ is smooth. By adapting arguments of standard Fourier theory for Banach space valued, smooth, compactly supported functions,  
                    we have the following estimate by adapting part (d) of Theorem 7.4 of \cite{Rudin}:
   \be \label{bdd} \| \widehat{\lambda H_1}(\underline{\mu}, \underline{\xi})\| \leq C_k (1+\sum_r |\mu_r|^2+\sum_{jl} | \xi_{jl} |^2)^{-k},\ee for every $k=0,1,2, \ldots$, 
   with some constants $C_k$. This shows that the right hand side of (\ref{eqn1}) converges absolutely. Let $Z$ denote the  right hand side of Equation (\ref{eqn1}).
                For any given $\gamma \in \hat{\clq}_{\pi(e)},$ applying the Fourier inversion formula to 
                $h:=\lambda (\gamma \circ H_1) \in C_c^\infty(\IR^{n+n^2})$: 
  $$ h(\underline{x},\underline{v})=\int \hat{h}(\underline{\mu},\underline{\xi}){\rm exp}(-2 \pi i ( \sum_k \mu_kx_k  +\sum_{jl} \xi_{jl} v_{jl})) d\underline{\mu} 
  d\underline{\xi}.$$  If $\gamma(\alpha(\chi)(\pi(e)))$ is zero then $\gamma(Z)=0$ and by definition of $\eta(F)$, $\gamma(\eta(F)(e))=0$ too.    
  On the other hand, if $\gamma(\alpha(\chi)(\pi(e)))$ is nonzero, putting $\underline{v}=\underline{\tau}(e)\equiv (( \gamma(T^U_{ij}(e) ))$, 
  we have, \bean \lefteqn{\gamma(\eta(F)(e))  =\gamma(\alpha(\chi)(\pi(e))) h(\underline{x}(e), \underline{\tau}(e))}\\
  &=& \gamma(\alpha(\chi)(\pi(e))) \int \hat{h}(\underline{\mu},\underline{\xi}){\rm exp}\left( -2 \pi i ( \sum_k \mu_kx_k(\pi(e))  +\sum_{jl} \xi_{jl} \gamma(T^U_{jl}(e))\right) d\underline{\mu}d\underline{\xi}=\gamma(Z),
  \eean   
     Note that we have used the estimate (\ref{bdd}) to justify interchanging  $\gamma$ with the integral sign,
    which proves (\ref{eqn1}). As $T^U_{ij}(e)$ and $x_k(\pi(e))$ are smooth functions of $e$, we can easily see the Fr\'echet continuity of $e \mapsto \eta(F)(e)$, using once again 
     the bound (\ref{bdd}). \qed

  We continue to denote by $\eta_U$ the unique $\ast$-homomorphic extension of $\eta_U$ on $C_c(E_U)$.

     \bthm
    \label{smooth}
    There exists a faithful smooth action $\eta$ of $\clq$ on $O(M)$
    satisfying $\eta(F)(e) \in \clq_{\pi(e)}$ for all $F\in C^{\infty}(O(M))$ and $e \in O(M)$.    \ethm
    {\it Proof}:\\
     Choose and fix a finite cover $U_1, \ldots, U_r$ of $M$ by coordinate neighbourhoods, 
  a $C^\infty$ partition of unity $\chi_i,i=1, \ldots, r$ subordinate to this cover and
  define $$ \eta(F)=\sum_{i=1}^r \eta_{U_i}((\chi_i \circ \pi)F),$$ for $F \in C(O(M))$. As each of the maps $\eta_{U_i}$ is Fr\'echet continuous on 
  $C_c^\infty(O(U_i))$,  $\eta$ too is Fr\'echet continuous.
  
   For $F,G \in C(O(M))$, writing $\hat{\chi}_i=\chi_i \circ \pi$, $\eta_i=\eta_{U_i}$,  we have $\eta(FG)=\sum_i \eta_i(\hat{\chi}_i FG)=\sum_{ij} \eta_i(\hat{\chi}_i \hat{\chi}_j FG)
   =\sum_{ij} \eta_i(\hat{\chi}_i F)\eta_j(\hat{\chi}_j G)$
   by (iii) of Lemma \ref{density}. But this is nothing but 
   \begin{displaymath}
    \left( \sum_i  \eta_i(\hat{\chi}_iF) \right)\left( \sum_j \eta_j (\hat{\chi}_j G) \right)=
   \eta(F) \eta(G).
   \end{displaymath}
   Similarly we can show 
    $\eta(\overline{F})=\eta(F)^*$. Thus, $\eta$ is a $\ast$-homomorphism. Being Fr\'echet continuous, it maps $C^\infty(O(M))$ to $C^\infty(O(M), \clq)$. 
    
    To complete the proof, consider the Fr\'echet dense $\ast$-subalgebra $\cla$ mentioned before, on which $\alpha$ is algebraic. Then 
     the algebra  $\clb$ generated by $\theta^i_\omega$, $1\leq i \leq n$, $\omega \in \Omega^1(\cla)$, is Fr\'echet dense in $C^\infty(O(M))$. 
     We claim  that for $\omega=f\omega^\prime_j$, 
 where $f \in C_c^\infty(U)$ for some coordinate neighbourhood $U$ and $(\omega^\prime_1, \ldots, \omega^\prime_n)$ are $U$-orthonormal one-forms, 
 $\eta(\theta^i_\omega)(e)=\lgl\lgl(\beta_{i} \ot
1_\clq),d\alpha_{(1)}(\omega)\rgl\rgl(\pi(e))$,  $e=(\pi(e), \beta_1, \ldots, \beta_n).$ This can be verified from the definition. However, such one-forms 
  comprise a Fr\'echet dense  subspace, hence we get the above  for all $\omega$. From this, it is clear that       
     $$\eta(\theta^{i}_{\omega})=\theta^{i}_{\omega_{(0)}} \ot \omega_{(1)}$$ for all $\omega \in \Omega^1(\cla), $ using the Sweedler-type  notation. 
      As $\alpha$ is algebraic on $\cla$, we have $\eta(\clb) \subseteq \clb \otimes \clq_0$.  It The co-associativity of 
       $\eta$ on $\clb$ (hence also  on 
      $C(O(M))$ by density ) follows from the co-associativity of $\alpha$ and moreover, we have $({\rm id} \ot \epsilon)\circ \alpha={\rm id}$ on $\cla$ which implies 
      $({\rm id} \ot \epsilon)\circ \eta={\rm id}$ on $\clb$.  That is,  $\eta$ on $\clb$ is a Hopf-algebraic co-action, 
      so in particular we have ${\rm Sp}\eta(\clb)(1 \otimes \clq_0)=\clb \ot \clq_0$,
       which is Fr\'echet dense in $C^\infty(O(M), \clq)$.  Hence $\eta$ is a smooth action in our sense. 
       To show faithfulness of $\eta$, consider $1\leq p\leq r$, $f\in C^{\infty}(M)$ and $\chi\in C^{\infty}_{c}(U_{p})$ 
       such that $0\leq \chi\leq 1$, $\chi\chi_{p}f=\chi_{p}f.$ From Lemma \ref{orthogonality} we get
       \begin{eqnarray*}
        \alpha(\chi_{p}f)(m)= \sum_{j=1}^{r}\eta((\chi_{p}\circ\pi)t^{U_{p}}_{1j})(e)\eta((\chi f\circ\pi)t^{U_{p}}_{1j})(e)
       \end{eqnarray*}
       for all $e\in O(M)$ with $\pi(e)=m$. Thus for all $1\leq p\leq r$ and $f\in C^{\infty}(M)$, $\alpha(\chi_{p}f)(m)$ is contained in the $C^\ast$ algebra generated by $\{\eta(F)(e)|F\in C^{\infty}(O(M)), e\in O(M)\}$. Hence the $C^{\ast}$-algebra also contains the $C^{\ast}$ closure of $\{\alpha(f)(m)|f\in C^{\infty}(M), m\in M\}$ which is $\clq$ by faithfulness of $\alpha$.
Finally, $\eta(F)(e)\in\clq_{\pi(e)}$ because $T_{ij}^{U}(e)\in \clq_{\pi(e)}$ by construction.
               \qed

\section{Isometric actions, i.e. actions commuting with the Laplacian}
\subsection{Definition and smoothness of isometric action}
We now discuss the apparently stronger condition of isometry, as in \cite{Goswami}. Throughout this section, let 
$M$ be a compact Riemannian manifold without boundary, with the Riemannian volume form $dvol$, $\clh:=L^2(M, dvol)$ and let $\tau$ denote the functional on $C(M)$ given by 
 $\tau(f)=\int_M f dvol$. 
 We denote by 
$\cll$   the Hodge Laplacian $-d^*d$ to $0$-forms, which is a self-adjoint operator on $\clh$. It is known ( see e.g. \cite{Donnedly} and the references therein) 
that $\cll$ has discrete spectrum given by 
 eigenvalues, say $\{ \lambda_i, i\geq 1 \}$ having finite multiplicities and the eigenvectors are in fact smooth functions. Let 
$\{e_{ij}:j=1,...,d_{i}\}$ be the orthonormal eigenvectors  of $\cll$ forming a
basis for
the eigen space corresponding to the eigenvector $\lambda_{i}$. We denote the linear span of 
$\{e_{ij}:1\leq j\leq d_{i},i\geq 1\}$ by $\cla_{0}^{\infty}$, which is a
subspace  of $C^{\infty}(M)$. Clearly, $\cll$ maps $C^\infty(M)$ to itself and we denote the restriction on $\cll$ to $C^\infty(M)$ (which is a Fr\'echet continuous operator) 
 again by the same symbol.
  
\bdfn
An action $\alpha$ of a CQG $\clq$ on $C(M)$ where $M$ is  a compact manifold $M$ without boundary, is said to be isometric if 
$\alpha(C^\infty(M)) \subseteq C^\infty(M, \clq)$ and for every state $\phi$ on $\clq$,
the map $({\rm id}\ot \phi)\alpha$ commutes with $\cll= -d^{*}d$ on $C^\infty(M)$. 
\edfn 
 We have
the following:
 \bthm
 \label{iso_smooth}
Any isometric action $\alpha$ is smooth, i.e. the linear span of $\alpha(C^\infty(M))(1 \ot \clq)$ is Fr\'echet dense in $C^\infty(M, \clq)$.
 \ethm
{\it Proof}:\\
Clearly,  $\alpha$ is
algebraic over $\cla_{0}^{\infty}$, hence $Sp \ \alpha(\cla_{0}^{\infty})(1\ot\clq_{0})=\cla_{0}^{\infty}\ot\clq_{0}$. It suffices to show that 
$\cla_{0}^{\infty}$ is Fr\'echet dense in $C^{\infty}(M)$. By Theorem 1.2 of \cite{Donnedly} there are constants $C$ and $C^{\prime}$
such that $||e_{ij}||_{\infty}\leq C|\lambda_{i}|^{\frac{n-1}{4}}$ and
$d_{i}\leq C^{\prime}|\lambda_{i}|^{\frac{n-1}{2}}$, where $n$ is the dimension of 
the manifold. For $f\in C^{\infty}(M)$ there are complex numbers $f_{ij}$ such
that $\sum_{ij} f_{ij}e_{ij}$ converges to $f$ in $L^{2}$ norm. Since
$f\in {\rm Dom}(\cll^{k})$ for all $k\geq 1$,
$\sum_{ij}|\lambda_{i}|^{2k}|f_{ij}|^{2}<\infty$ for all $k$. Choose and fix
sufficiently large $k$ such that $\sum_{i\geq 0}|\lambda_{i}|^{n-2k}<\infty$. This is possible by the well-known Weyl asymptotics of the eigenvalues of Laplacian.

Now, $\cll(\sum_{i \leq N, j \leq d_i}f_{ij}e_{ij})=\sum_{i\leq N, j \leq d_i}\lambda_{i}f_{ij}e_{ij}$
converges to $\cll(f)$ in the $L^{2}$ norm as $N \raro \infty$. By the Cauchy-Schwartz inequality,
$$\sum_{ij}|\lambda_{i}f_{ij}|||e_{ij}||_{\infty}\leq
C(C^{'})^{\frac{1}{2}}(\sum_{ij}|f_{ij}|^{2}|\lambda_{i}|^{2k})^{\frac{1}{2}}
(\sum_{i}|\lambda_{i}|^{n-2k})<\infty.$$
Hence $\lim_{N\raro \infty} \| \cll(\sum_{i \leq N, j \leq d_i}f_{ij}e_{ij})-\cll(f)\|_\infty =0$.  Similarly we can show that
$$\lim_{N \raro \infty} \cll^{k}(\sum_{i \leq N, j \leq d_i}f_{ij}e_{ij})=\cll^k(f)$$   in the sup norm of $C(M)$ for any
$k \geq 1$. Hence  $\cla_{0}^{\infty}$ is Fr\'echet dense in $C^{\infty}(M)$.
 \qed \\
 We also have the following:
 \blmma
  Any isometric action preserves the corresponding Riemannian inner product as well as the Riemannian volume measure. 
  \elmma
{\it Proof:}\\
 It is enough to 
 prove the following: \be \label{753} \lgl\lgl da_{(0)},db_{(0)}\rgl\rgl\ot  a_{(1)}^*b_{(1)}=\alpha(\lgl\lgl da,db\rgl\rgl),\ee
 where $\lgl\lgl df,dg\rgl\rgl=\cll(\overline{f}g)-\cll(\overline{f})g-\overline{f}\cll(g)$, $\forall f,g$ belonging to the subalgebra  $\cla$  as in Proposition \ref{smooth_action_dense}.
 However, it is straightforward to verify (\ref{753}) using $\cll(f_{(0)})\ot f_{(1)}=\alpha(\cll(f))$ for $f \in \cla$ as well as the fact that $\alpha$ is a $\ast$-homomorphism.
 
 The proof of the statement about Riemannian volume measure preservation can be found in Lemma 2.5 of \cite{Goswami}. \qed

  \subsection{ Commutativity of higher order partial derivatives}

 Let $\nabla$ be the Levi-Civita connection viewed as a map from $\Omega^{1}(C^{\infty}(M))$ to $\Omega^{1}(C^{\infty}(M))\ot_{C^{\infty}(M)} \Omega^{1}(C^{\infty}(M))$.  
 We have the following, using the observations that $\omega(X_{dh})=\lgl\lgl\omega,dh\rgl\rgl$, $df(Z)=\lgl\lgl X_{df},Z \rgl\rgl$ and $\lgl\lgl\nabla(df),dg\ot dh\rgl\rgl=\lgl\lgl\nabla_{X_{dh}}(df),dg\rgl\rgl$:

  \begin{eqnarray} \label{levi_civita_formula}
  \lgl\lgl\nabla(df), dg \otimes_{C^{\infty}(M)} dh\rgl\rgl&=&\frac{1}{2}( -\lgl\lgl df, d(\lgl\lgl dg, dh\rgl\rgl)\rgl\rgl+\lgl\lgl dh, d(\lgl\lgl df,dg\rgl\rgl)\rgl\rgl \nonumber \\
	&+&\lgl\lgl dg, d(\lgl\lgl df,dh\rgl\rgl)\rgl\rgl ),
	\end{eqnarray}
  where $f,g,h\in C^{\infty}(M)_{\mathbb{R}}$. 
  This can be derived by a slightly long but straightforward calculation using the standard formula for the Levi-Civita connection on vector fields (see, for example, page 69 of 
   \cite{Lee}), (\ref{x_omega}) as well as the formulae $\nabla_X(\omega)(Y):=X\omega(Y)-\omega(\nabla_X(Y))$ and $d\omega(X,Y)=X\omega(Y)-Y\omega(X)-\omega([X,Y])$ for all $X,Y\in\chi(M)$ (see page 54 of \cite{Lee}).


 \blmma
 \label{cov_2}
 Let $f,g,h\in\cla_{\rm s.a.}$. Then \begin{eqnarray} &&\lefteqn{\lgl\lgl dg_{(0)},d(\lgl\lgl df_{(0)},dh_{(0)}\rgl\rgl)\rgl\rgl\ot g_{(1)}f_{(1)}h_{(1)}}{\nonumber}\\
 &=&\lgl\lgl dg_{(0)},d(\lgl\lgl df_{(0)},dh_{(0)}\rgl\rgl)\rgl\rgl\ot f_{(1)}g_{(1)}h_{(1)}\label{030}\\
 &=& \lgl\lgl dg_{(0)},d(\lgl\lgl df_{(0)},dh_{(0)}\rgl\rgl)\rgl\rgl\ot h_{(1)}g_{(1)}f_{(1)}\label{040}
 \end{eqnarray}
 \elmma
 {\it Proof}:\\
We have
\begin{eqnarray*} 
&&\lgl\lgl dg_{(0)},d(\lgl\lgl df_{(0)},dh_{(0)}\rgl\rgl)\rgl\rgl\ot g_{(1)}f_{(1)}h_{(1)}\\
&=& \left(\cll(g_{(0)}\lgl\lgl df_{(0)},dh_{(0)}\rgl\rgl)-\cll(g_{(0)})\lgl\lgl df_{(0)},dh_{(0)}\rgl\rgl-g_{(0)}\cll(\lgl\lgl df_{(0)},dh_{(0)}\rgl\rgl)\right)\ot g_{(1)}f_{(1)}h_{(1)} \end{eqnarray*}
We compute the terms individually, using ${\mathcal L}(\phi_{(0)}) \otimes \phi_{(1)}
        =\alpha({\mathcal L}(\phi)).$
\begin{eqnarray*}
&& \cll(g_{(0)}\lgl\lgl df_{(0)},dh_{(0)}\rgl\rgl)\ot g_{(1)}f_{(1)}h_{(1)}\\
&=& (\cll\ot {\rm id})(\alpha(g)\alpha(\lgl\lgl df,dh\rgl\rgl))\\
&=& (\cll\ot {\rm id})\alpha(\lgl\lgl df,gdh\rgl\rgl)\\
&=& \cll(g_{(0)}\lgl\lgl df_{(0)},dh_{(0)}\rgl\rgl)\ot f_{(1)}g_{(1)}h_{(1)}.
\end{eqnarray*}
\begin{eqnarray*}
&&\cll(g_{(0)})\lgl\lgl df_{(0)},dh_{(0)}\rgl\rgl \ot g_{(1)}f_{(1)}h_{(1)}\\
&=& (\cll\ot{\rm id})(\alpha(g))\alpha(\lgl\lgl df,dh\rgl\rgl )\\
&=& \alpha(\cll(g))\alpha(\lgl\lgl df,dh\rgl\rgl )\\
&=& \alpha(\lgl\lgl df\cll(g),dh\rgl\rgl )\\
&=& \lgl\lgl df_{(0)}\cll(g_{(0)}),dh_{(0)}\rgl\rgl \ot f_{(1)}g_{(1)}h_{(1)}.
\end{eqnarray*}
Recall that $\tau(f)=\int f {\rm dvol}$. As  $\cll$ is a self adjoint operator on $\clh$, for  $\phi\in C^{\infty}(M)$ we get
\begin{eqnarray*}
&&\lgl\lgl\cll(\lgl\lgl df_{(0)},dh_{(0)}\rgl\rgl )g_{(0)}\ot g_{(1)}f_{(1)}h_{(1)},\alpha(\phi)\rgl\rgl\\
&=& \tau\left(\cll(\lgl\lgl df_{(0)},dh_{(0)}\rgl\rgl) g_{(0)}\phi_{(0)}\right)h_{(1)}f_{(1)}g_{(1)}\phi_{(1)}\\
&=& \tau\left(\lgl\lgl df_{(0)},dh_{(0)}\rgl\rgl \cll(g_{(0)}\phi_{(0)})\right)h_{(1)}f_{(1)}g_{(1)}\phi_{(1)} \ (by \ self-adjointness \ of \ \cll)\\
&=& \tau\left(\lgl\lgl df_{(0)},dh_{(0)}\rgl\rgl (\cll(g_{(0)})\phi_{(0)}+g_{(0)}\cll(\phi_{(0)})+\lgl\lgl dg_{(0)},d\phi_{(0)}\rgl\rgl)\right)h_{(1)}f_{(1)}g_{(1)}\phi_{(1)}.
\end{eqnarray*}
Using $\alpha(\lgl\lgl df,dh\rgl\rgl)=\lgl\lgl df_{(0)},dh_{(0)}\rgl\rgl\ot f_{(1)}h_{(1)}=\lgl\lgl df_{(0)},dh_{(0)}\rgl\rgl\ot h_{(1)}f_{(1)}$ as observed earlier, we get:
\begin{eqnarray*}
&&\tau\left(\lgl\lgl df_{(0)},dh_{(0)}\rgl\rgl\cll(g_{(0)})\phi_{(0)}\right)h_{(1)}f_{(1)}g_{(1)}\phi_{(1)}\\
&=& (\tau\ot{\rm id})\left(\alpha(\lgl\lgl df,dh\rgl\rgl\cll(g)\phi)\right)\\
&=& (\tau\ot{\rm id})\left(\alpha(\lgl\lgl dh,df\rgl\rgl\cll(g)\phi)\right)\\
&=& (\tau\ot{\rm id})\left(\alpha(\lgl\lgl dh,\cll(g)df\rgl\rgl\phi)\right)\\
&=& \tau\left(\lgl\lgl dh_{(0)},df_{(0)}\rgl\rgl \cll(g_{(0)})\phi_{(0)}\right)\ot h_{(1)}g_{(1)}f_{(1)}\phi_{(1)}.
\end{eqnarray*}
Similarly, \begin{displaymath}\tau\left(\lgl\lgl df_{(0)},dh_{(0)}\rgl\rgl g_{(0)}\cll(\phi_{(0)}\right)h_{(1)}f_{(1)}g_{(1)}\phi_{(1)}=\tau(\lgl\lgl df_{(0)},dh_{(0)}\rgl\rgl g_{(0)}\cll(\phi_{(0)})\ot h_{(1)}g_{(1)}f_{(1)}\phi_{(1)}.\end{displaymath}
Also,
\begin{eqnarray*}
&& \tau\left(\lgl\lgl df_{(0)},dh_{(0)}\rgl\rgl\lgl\lgl dg_{(0)},d\phi_{(0)}\rgl\rgl\right)h_{(1)}f_{(1)}g_{(1)}\phi_{(1)}\\
&=& (\tau\ot {\rm id})\left( \lgl\lgl(d\ot {\rm id})\alpha(h),(d\ot {\rm id})\alpha(f)\rgl\rgl\lgl\lgl (d\ot {\rm id})\alpha(g),(d\ot {\rm id})\alpha(\phi)\rgl\rgl\right).
\end{eqnarray*}
By (\ref{comm111}), 
\begin{eqnarray*}
&&\tau\left( \lgl\lgl df_{(0)},dh_{(0)}\rgl\rgl \lgl\lgl dg_{(0)},d\phi_{(0)}\rgl\rgl\right)h_{(1)}f_{(1)}g_{(1)}\phi_{(1)}\\
&=& \tau \left(\lgl\lgl df_{(0)},dh_{(0)}\rgl\rgl\lgl\lgl dg_{(0)},d\phi_{(0)}\rgl\rgl\right)h_{(1)}g_{(1)}f_{(1)}\phi_{(1)},
\end{eqnarray*}
hence \begin{displaymath}\lgl\lgl g_{(0)}\cll(\lgl\lgl df_{(0)},dh_{(0)}\rgl\rgl)\ot g_{(1)}f_{(1)}h_{(1)},\alpha(\phi)\rgl\rgl
=\lgl\lgl g_{(0)}\cll(\lgl\lgl df_{(0)},dh_{(0)}\rgl\rgl)\ot f_{(1)}g_{(1)}h_{(1)},\alpha(\phi)\rgl\rgl\end{displaymath} for all $\phi\in\cla$. As ${\rm Sp} \ \{ \alpha(\phi)q:~\phi \in \cla, q \in \clq\}$
 is dense in the Hilbert module $\clh\bar{\ot}\clq$, we get
$$g_{(0)}\cll(\lgl\lgl df_{(0)},dh_{(0)}\rgl\rgl)\ot g_{(1)}f_{(1)}h_{(1)}=g_{(0)}\cll(\lgl\lgl df_{(0)},dh_{(0)}\rgl\rgl)\ot f_{(1)}g_{(1)}h_{(1)}.$$
Combining all these we get (\ref{030}).

To prove the other equality,  note that  $$\lgl\lgl df_{(0)},dh_{(0)}\rgl\rgl \ot f_{(1)}h_{(1)}=\lgl\lgl df_{(0)},dh_{(0)}\rgl\rgl \ot h_{(1)}f_{(1)}=\lgl\lgl dh_{(0)},df_{(0)}\rgl\rgl \ot h_{(1)}f_{(1)}.$$
 Hence \bean 
 &&\lefteqn{\lgl \lgl  dg_{(0)}, d\lgl\lgl df_{(0)},dh_{(0)}\rgl\rgl \rgl \rgl \ot  g_{(1)}f_{(1)}h_{(1)}}\\
 &=&\lgl \lgl  dg_{(0)},d \lgl\lgl dh_{(0)},df_{(0)}\rgl\rgl \rgl \rgl\ot g_{(1)} h_{(1)}f_{(1)}\\
  &=& \lgl \lgl  dg_{(0)},d \lgl\lgl dh_{(0)},df_{(0)}\rgl\rgl \rgl \rgl\ot h_{(1)} g_{(1)}f_{(1)},
 \eean
 where in the last step we have used (\ref{030}) interchanging $f$ and $h$.
  \qed 

We consider $(\nabla \overline{\ot} {\rm id}_\clq): \Omega^1(C^\infty(M)) \overline{\ot} \clq \raro \Omega^2(C^\infty(M)) \overline{\ot} \clq$ 
as follows. Let $m \in M$, $(U, x_1, \ldots, x_n)$ a local chart around $m$ and $\Omega \in \Omega^1(C^\infty(M)) \overline{\ot} \clq$ such that 
$\Omega(x)=\sum_i dx_i|_x \Omega_i(x)$ $\forall x \in U$. Then we define 
\begin{displaymath}(\nabla \overline{\ot} {\rm id}_\clq)(\Omega)(m)=\sum_i \left( (\nabla(dx_i)(m) \ot 1) \Omega_i(m)+dx_i(m) \ot (d\Omega_i)(m) \right).\end{displaymath}
By standard arguments one can see that the above does not depend on the choice of local coordinates and $(\nabla \overline{\ot} {\rm id}_\clq) $ is Fr\'echet continuous.

\bcrlre
\label{connection_pres_lemma}
For $\omega \in \Omega^1(C^\infty(M)),$ we have 
\be \label{connection_pres_2} (\nabla \overline{\otimes} {\rm id} )(d\alpha_{(1)}(\omega)) =U^{(2)}(\nabla(\omega)),\ee
where $U^{(2)}$ is the equivariant unitary representation on $\Omega^2(C^\infty(M))$ constructed  in Lemma  \ref{lift_eq_rep}.
\ecrlre
{\it Proof:}\\
It suffices to prove \be \label{prelim} (\nabla \ot {\rm id})(d\alpha_{(1)}(dg))=U^{(2)}(\nabla(dg)) \ee for $g \in \cla_{\rm s.a.}$. Indeed, we can prove the  identity
 (\ref{connection_pres_2}) from 
(\ref{prelim}) in two steps. First, we prove it for $\omega=(dg) f$, where $f,g \in \cla_{\rm s.a.}$, using 
the Leibniz rule for connection and the fact $U^{(2)}(\Theta f)=U^{(2)}(\Theta) \alpha(f)$ $\forall \Theta \in \Omega^2(C^\infty(M))$. Then 
 we use the Fr\'echet-continuity of $\nabla \overline{\ot} {\rm id}$ and 
 Fr\'echet-density of the complex linear span of one-forms of the form $(dg)f,$ $f,g \in \cla_{\rm s.a.}$ in $\Omega^1(C^\infty(M))$.

As the complex linear span of elements of the form $(df \otimes_{\cla} dh) \phi,$ where $f,h,\phi \in \cla_{\rm s.a.}$ is dense in $\Omega^2(C^\infty(M))$, 
the complex linear span  of $ U^{(2)}((df \otimes_{\cla} dh) \phi)q=df_{(0)} \ot_{\cla} dh_{(0)} \phi_{(0)} \ot f_{(1)}h_{(1)}\phi_{(1)}q $, $q \in \clq$, 
is dense in $\Omega^2(C^\infty(M)) \overline{\ot} \clq$, i.e. elements of the form $U^{(2)}(df \ot_{\cla} dh)$ constitute a right $C^\infty(M, \clq)$-total subset. 
To prove (\ref{prelim}) it is enough to show that $\lgl\lgl(\nabla (dg_{(0)}) \ot dg_{(1)}, \Omega\rgl\rgl=\lgl\lgl U^{(2)}(\nabla(dg)), \Omega\rgl\rgl$ for
all $\Omega $ in some subset whose right $C^\infty(M, \clq)$-linear span is dense in $\Omega^2(C^\infty(M)) \overline{\ot} \clq$, in particular, 
  \begin{displaymath}\Omega=U^{(2)}(df \ot_{\cla} dh)=df_{(0)} \ot_{\cla} dh_{(0)} \ot f_{(1)}h_{(1)}.\end{displaymath} But we have the following by 
  combining the formula (\ref{levi_civita_formula}), Lemma \ref{cov_2} and the equivariance of $U$:

\be \label{connection_pres_1}\lgl\lgl\nabla(dg_{(0)}), df_{(0)} \otimes_{\cla} dh_{(0)} \rgl\rgl\ot g_{(1)}f_{(1)}h_{(1)}=\alpha(\lgl\lgl\nabla(dg),df \otimes_{\cla} dh \rgl\rgl). \ee
Indeed, as $U^{(2)}$ is equivariant, the right hand side of the above is equal to \begin{displaymath}\lgl\lgl U^{(2)}(\nabla(dg)), U^{(2)}(df \ot_{\cla} dh)\rgl\rgl,\end{displaymath} which completes the proof. \qed

   This leads to the following:
   \bthm 
   \label{commutativity_all}
   For any $m \in M$, local coordinates $(W,(x_1, \ldots, x_n))$ around $m$, a positive integer $k$  and $f \in C^\infty(M)$, we have 
   $$(\frac{\partial}{\partial x_{i_1}} \ldots \frac{\partial}{\partial x_{i_k}} \alpha(f))(m) \in \clq_m$$
    $\forall i_l \in \{ 1, \ldots, n \}.$
      \ethm
   {\it Proof:}\\ Let $\cla\subset C^{\infty}(M)$ be as in Proposition \ref{smooth_action_dense} and $\omega_1, \ldots, \omega_n$ be $W$-orthonormal smooth one forms 
   where $W$ is some coordinate neighbourhood 
    around $m$. Write $X_i=X_{\omega_i}$, using the notation introduced in Subsection \ref{form}.
   The result is clearly equivalent to the following : $ (X_{i_1} \ldots X_{i_k} (\alpha(f)))(m) \in \clq_m$ for $k \geq 1$, $1 \leq i_j \leq  n$.
  
   Let us define the following maps $\nabla^k$ from $ \Omega^1(C^\infty(M)) $ to $\Omega^{k+1}(C^\infty(M))$ for $k=1,2,\ldots.$  Let $\sigma_{ij}$
    denote the map which flips the $i$-th and the $j$-th copies of $\Omega^1(C^\infty(M))$. Define  $\nabla^1=\nabla.$
    Then, consider the $\IC$-linear map $ T:\Omega^1(C^\infty(M)) \ot_{\IC} \Omega^1(C^\infty(M))\raro \Omega^3(C^\infty(M))$ defined by 
     \begin{displaymath}T(\omega \ot_{\IC} \eta)=\sigma_{23} (\nabla(\omega) \ot_{C^\infty(M)}\eta)+\omega \ot_{C^\infty(M)} \nabla(\eta).\end{displaymath} 
    Verify using the Leibniz rule for connections that $T(\omega f \ot_{\IC} \eta)=T(\omega \ot_{\IC} f \eta)$, hence it descends to a map
      say $\overline{T}$ on $\Omega^2(C^{\infty}(M))$. We define $\nabla^{2}=\overline{T}\circ\nabla$. In  a similar way, for $k \geq 2$, define $\nabla^{k}:=(\sum_{l=1}^{k-1} \sigma_{l+1 k+1} \circ \nabla_l +
    \nabla_k) \circ \nabla^{k-1}$, where $\nabla_l$ is the map which acts on the $l$-th copy of $\Omega^1(C^\infty(M))$ by $\nabla$ leaving the other copies unaffected.

    It follows by repeated application of (\ref{connection_pres_2}) and the commutativity of $\clq_m$ that 
		\begin{displaymath}\nabla^k(df_{(0)})(m) \ot f_{(1)} =
    U^{(k+1)}(\nabla^k(df))(m)\in (T^*_m M)^{\ot^{k+1}} \ot \clq_m.\end{displaymath} On $W$, 
     we can write $df_{(0)} \ot f_{(1)}=\sum_{i=1}^n \omega_i X_i(f_{(0)}) \ot f_{(1)}$, hence \begin{displaymath}\nabla(df_{(0)}) \ot f_{(1)}=\sum_i  \nabla(\omega_i) X_i(f_{(0)}) \ot 
    f_{(1)}+\sum_{i,j=1}^n \omega_i \ot_\cla \omega_j X_jX_i(f_{(0)}) \ot f_{(1)}.\end{displaymath} Evaluating at $m$ and noting that the first term belongs to $T^*_m M \ot T^*_m M \ot \clq_m$ 
     by definition of $\clq_m$, we conclude that the second term must belong to this space too. Thus, taking inner product with $\omega_i(m) \ot \omega_j(m)$ for 
      any fixed $i,j$, we get $X_iX_j (\alpha(f))(m) \in \clq_m$. Expanding $\nabla^2$ in a similar way and using $X_i \alpha(f)(x), X_i X_j \alpha(f) (x)$ 
      are in $\clq_x$ for all $x \in W$, 
       we show $X_iX_jX_r \alpha(f)(m) \in \clq_m$ for all $i,j,r$. Proceeding inductively, using the expansion of $\nabla^k(df_{(0)}) \ot f_{(1)}$, 
       the statement of the theorem follows for any $k \geq 1$. \qed

\bcrlre
\label{higher_comm_cor}
For any $\omega, \eta \in \Omega^1(M)$, we have $(\frac{\partial}{\partial x_i}\lgl\lgl \eta \ot 1, d \alpha (\omega)\rgl\rgl)(m) \in \clq_m$.
\ecrlre
{\it Proof :} It is enough to prove it for $\omega=fdg,$ $f,g \in \cla$. Indeed 
\begin{displaymath}(\frac{\partial}{\partial x_i}\lgl\lgl \eta \ot 1, d \alpha (f dg)\rgl\rgl)(m)=\frac{\partial}{\partial x_i}\left(\alpha(f)X_\eta(\alpha(g))\right)(m) \in \clq_m,\end{displaymath} 
as $\frac{\partial}{\partial x_i} (X \alpha(g)) (m)\in \clq_m$ for every smooth vector field $X$ by Theorem \ref{commutativity_all}. \qed

\section{Main result}
   \blmma
  \label{affine}
 Let $\Phi$ be a smooth action of a CQG $\clq$ on a compact connected subset $W$ which is the closure of a bounded smooth domain, i.e. bounded, open connected subset of $R^N$ with smooth 
            boundary.  Suppose 
  furthermore that the action preserves  the usual (Euclidean) Riemannian inner product and for any $y \in W$, the algebra generated by 
   $\{ \Phi(g)(y), (\frac{\partial}{\partial y_{i_1}} \ldots \frac{\partial}{\partial y_{i_k}} \Phi(f))(y):~k \geq 1,1 \leq i_j \leq N,~f,g \in C^\infty(W)\}$
    is commutative, where  $y_1, \ldots, y_N$ denote the standard coordinates of $\IR^N$.  Then $\Phi$ is
affine
i.e.
\begin{eqnarray}
\label{exp_phi}
 \Phi(y_{i})=1\ot q_{i}+ \sum_{j=1}^{N}y_{j}\ot q_{ij}, \ for \ some \
q_{ij},q_{i}\in \clq,
\end{eqnarray}
  for all $i=1,...,N$. 
  \elmma
 {\it Proof:}\\
 Let $W=\overline{V}$, where $V$ is a bounded, open connected set with smooth boundary. Note that $(dy_1, \ldots, dy_N)$ is an orthonormal basis of $T^*_mW$ at every point $m $. Let $D^k_i(m)=\frac{\partial}{\partial y_i}|_m \Phi(y_k)$, 
   $D^k_{ij}(m)=\frac{\partial^{2}}{\partial y_i\partial y_j}|_m \Phi(y_k)$. As $V$ is open and connected, it suffices to 
    prove $D^k_{ij}(m)=0$ for all $m \in V$. Using the arguments and discussion in the beginning of Subsection 3.2 (page 10), we get a $\Phi$-equivariant unitary representation $\Gamma:=d\Phi_{(1)} $ as in the Definition 3.6,
    satisfying $\Gamma(df) =d \Phi(f)$ for every $f \in C^\infty(W)$, in particular, $\Gamma(dy_i)=\sum_{j=1}^N dy_j D^i_j.$ As the Hilbert $C^\infty(W)$-module of one-forms of $W$ is 
  free of rank $N$ with the (orthonormal w.r.t. the Euclidean Riemmanian structure)  basis $dy_1, \ldots, dy_N$,  we can adapt the arguments of of Lemma \ref{unitary} 
  to conclude that 
  $(( D^i_j ))_{i,j=1}^N $ is a unitary element of $M_N(C^\infty(W))$. From the unitarity as well as self-adjointness of $D^i_j$'s (which follows because $y_i$'s are self-adjoint), we 
   get the following two equations:
    \be \label{(20)}  \sum_{l=1}^N D^l_i D^l_j=\delta_{ij}1_{\mathcal Q}, \ee
    \be \sum_{l=1}^N D^i_l D^j_l =\delta_{ij}1_{\mathcal Q}.\ee
   
 Applying $\frac{\partial}{\partial y_{k}}$ to equation (\ref{(20)}), and using the commutativity of  $D^l_{jk}$ and $D^l_i$  for all choices of the upper and lower indices, we get
 \begin{eqnarray}
 \label{(21)}
 \sum_{l=1}^{N}(D^l_{ik}D^l_{j}+D^{l}_{jk}D^{l}_{i})=0.
 \end{eqnarray}
 Now we denote the $N^{2}\times N$ matrix whose $(ij)k$-th entry $A_{(ij),k}$ is $D^{k}_{ij}(y)$ by $A$ and $N\times N$ matrix whose $ij$-th entry $B_{ij}$ is
  given by $D^{i}_{j}(y)$ by $B$. Then the $(ij)k$-th entry of the matrix $C=AB$ denoted by $C_{(ij)k}$ is given by $\sum_{l=1}^{N}A_{(ij),l}D^{l}_{k}$.  From the 
  equation (\ref{(21)}), we get for all $i,j,k$,
 \begin{eqnarray}
 \label{(22)}
 C_{(ik)j}+C_{(jk)i}=0.
 \end{eqnarray}
 Now we observe that $C_{(ij)k}=C_{(ji)k}$ for all $i,j,k$. Hence by repeated application of equation (\ref{(22)}),
 \begin{displaymath}
 C_{(ik)j}= C_{(ki)j}=-C_{(ji)k}=-C_{(ij)k}=C_{(kj)i}=C_{(jk)i}.
 \end{displaymath}
 Again by equation (\ref{(22)}), we get $C_{(ik)j}=0$ for all $i,j,k$ i.e. $C=0$. As $B$ is unitary, we conclude that $A=0$, i.e. $D^{k}_{ij}$ is zero for all $i,j,k$ completing the proof of the Lemma. \qed
 \bcrlre
 \label{split_main}
 In the set up of Lemma \ref{affine}, the $C^{\ast}$ algebra generated by $\{q_{i},q_{kl}:i,k,l=1,...,n\}$ is commutative. In particular, if $\Phi$ is faithful, $\clq$ must be 
  commutative.
 \ecrlre
 {\it Proof}:\\
 Observe that $(\frac{\partial}{\partial y_{j}}\Phi(y_{i}))(x)$ and $\Phi(y_{k})(x)$ must commute for $i,j,k=1,...,n$ and for all $x\in W$ as $\Phi$ is isometric. But we have $(\frac{\partial}{\partial y_{j}}\Phi(y_{i}))=1\ot q_{ij}$. Thus it follows from the expression of $\Phi(y_{k})$ given by (\ref{exp_phi}), 
 \begin{eqnarray}
 \label{quad}
 1\ot q_{k}q_{ij}+\sum_{l=1}^{N}y_{l}\ot q_{kl}q_{ij}=1\ot q_{ij}q_{k}+\sum_{l=1}^{N}y_{l}\ot q_{ij}q_{kl}.
 \end{eqnarray}
 Clearly the set $\{1,y_{1},...,y_{N}\}$ is a linearly independent set as $W$ has a non-empty interior, hence we get the commutativity of $q_{k}, q_{ij}$ and $q_{kl},q_{ij}'s$ by comparing the coefficients of $\{1,y_{1},...,y_{N}\}$ on both sides of \ref{quad}.\qed
 
\indent Now, we are in a position to prove the main theorem of the paper. 

\bthm 
\label{main_thm}
Let $\alpha$ be  a faithful isometric action of a CQG $\clq$ on a compact, connected, Riemannian $n$-manifold $M$. Then $\clq$ is commutative. 
 \ethm
{\it Proof:}\\ 
We break the proof into several steps. As the CQG $\clq$ is commutative if and only if the corresponding reduced CQG $\clq_{r}$ is so, we may assume without loss of generality 
that $\clq=\clq_{r}$ is reduced CQG so that Theorem \ref{automatic_lift} applies. \\
{\it Step 1}\\
Let $P=O(M)$ and let $\eta$ be the lift of the action on $P$ obtained by Theorem \ref{smooth}, which 
satisfies $\eta(\phi)(e) \in \clq_{\pi(e)}$ $\forall \phi \in C^\infty(P)$ and $e \in P$.
 We claim that 
$X (\eta(\phi))(e) \in \clq_{\pi(e)}$ too for all $X\in\chi(P)$. Once this is proved, we can apply Theorem \ref{automatic_lift} (as $\clq_{\pi(e)}$ is commutative) 
to get  some Riemannian structure on $P$ for which $\eta$ is inner product preserving. To prove the claim, it is enough to consider $\phi$ of the form 
$(f \circ \pi)t_{ij}^{(U,\underline{\omega}^\prime)}$ in the notation of Subsection 3.3. Moreover, as in the proof of Lemma \ref{density}, let us choose some local coordinate 
$(V, (x_1, \ldots, x_n))$ for $M$ around $\pi(e)$,  $V$-orthonormal one-forms $\omega^{\prime \prime}_1, \ldots, \omega^{\prime \prime}_n$ and the corresponding 
embedding of $\pi^{-1}(V)$ into $\IR^{n+n^2}$ using $x_r, t^{(V, \underline{\omega^{\prime \prime}})}_{pq},~r,p,q=1, \ldots, n$. 
 Denoting the canonical  coordinate functions 
  of $\IR^{n^2}$ by $y_{pq}, p,q =1,\ldots, n$, it suffices to verify the claim for $X$ of the form $\partial_r \equiv \frac{\partial}{\partial x_r}$ and $\partial_{pq}=
  \frac{\partial}{\partial y_{pq}}.$  But $$ \eta(\phi)(e)=\alpha(f)( \pi(e)) T^{(U,\underline{\omega^\prime})}_{ij}(e)
=\sum_{k=1}^{n} \alpha(f)(\pi(e)) t^{(V, \underline{\omega^{\prime \prime}})}_{ik}(e)\lgl\lgl \omega^{\prime \prime}_{k} \ot 1, d\alpha_{(1)}(\omega^{\prime}_{j})\rgl\rgl(\pi(e)).$$
From this expression, it is clear that $\partial_{pq}\eta(\phi)(e)$ is a complex linear combination of 
$$\{\alpha(f)(\pi(e)),~ \lgl\lgl \omega^{\prime \prime}_k \ot 1, d\alpha_{(1)}(\omega^\prime_j)\rgl\rgl(\pi(e)),~k,j=1, \ldots, n\} \subseteq  \clq_{\pi(e)}.$$ On the other hand, $\partial_r \eta(\phi)(e)$ is a complex linear 
 combination of $\partial_r(\alpha(f))(\pi(e)) (\in \clq_{\pi(e)})$ as well as 
 $\partial_r(\lgl\lgl \omega^{\prime \prime}_k \ot 1, d\alpha_{(1)}(\omega^\prime_j)\rgl\rgl)(\pi(e))$ ($k,j=1, \ldots n)$  
 which also belong to $\clq_{\pi(e)}$ by Corollary \ref{higher_comm_cor}.
\vspace{2mm}\\
{\it Step 2}:\\
We now fix a Riemannian structure on $P$ obtained by {\it Step 1} and want to lift $\eta$ further to a tubular normal neighbourhood of $P$. Consider a Fr\'echet
dense
subalgebra $\cld$ of $C^\infty(P)$ over which $\eta$ is
algebraic and $Sp\left(\eta(\cld)(1\ot\clq_0)\right)=\cld \ot\clq_0.$

As $P$ is parallelizable,  there is a Riemannian embedding (for the Riemannian structure of $P$ discussed in {\it Step 1}) into some $\IR^N$, so that 
 the corresponding  normal bundle is trivial. Recall from Lemma
\ref{global_diff} the global diffeomorphism $F$ and the corresponding isomorphism 
$$\pi_F:C^\infty(\cln_\epsilon P)\rightarrow C^\infty(P\times
B^{N-r}_\epsilon(0)),$$ where $r$ denotes the dimension of $P$. 
 Define 
$\widehat{\eta}: C^\infty( P\times
B^{N-r}_\epsilon(0)) \raro C^\infty (P\times
B^{N-r}_\epsilon(0), \clq)$ by $$ \widehat{\eta}(G)(e,b)=\eta(G_b)(e),$$ where $G_b: P \raro \IC$ is given by $G_b(e)=G(e,b)$. It is clearly a $C^{\ast}$ action 
and also satisfies
$\widehat{\eta}=\sigma_{23}\circ (\eta \ot {\rm id})$ on $\cld \ot C^\infty(B^{N-r}_\epsilon(0) )$, from which it is clear 
 that $\widehat{\eta}(\cld \ot C^\infty(B^{N-r}_\epsilon(0) ))(1 \ot \clq_0)=\cld \ot C^\infty(B^{N-r}_\epsilon(0) ) \ot \clq_0$. Thus, 
  $\widehat{\eta}$ is indeed a smooth action. 
Now we have 
$\pi_{F^{-1}}:C^\infty(M \times B^{N-r}_\epsilon(0))\rightarrow C^\infty(\cln_\epsilon P).$ Define  
$$\Phi:=(\pi_{F^{-1}}\ot {\rm id})\circ\widehat{\eta}\circ\pi_F:
C^\infty(\cln_\epsilon P)\rightarrow C^\infty(\cln_\epsilon
P,\clq).$$ We claim that 
$\Phi$ is a smooth action of $\clq$ on
$\cln_{\epsilon}P$.

 To see this, let $\widetilde{\cld}:=\pi_{F^{-1}}(\cld \ot C^\infty(B^{N-r}_\epsilon(0))),$ which is a 
Fr\'echet dense subalgebra of $C^\infty(\cln_{\epsilon}P).$ By construction,
$\Phi$ is algebraic over $\widetilde{\cld}$ and moreover,
 $Sp~\Phi(\widetilde{\cld})(1\ot\clq_0)=\widetilde{\cld}\ot\clq_0$ Also $\Phi$ is Fr\'echet continuous, hence smooth.\\
 {\it Step 3}:\\
 We claim that $\Phi$ preserves the Riemannian inner product of $\cln_\epsilon P.$ Denote the projection from $\cln_\epsilon P$ to $P$ by $\rho$.  As 
 the normal bundle is trivial, choose a smoothly varying
basis for normal space at each point of $P$. Let $y\in
\cln_{\epsilon}P$ and $\{e_{i}(y): i=1,\ldots,(N-r)\}$ be an orthonormal basis
for the normal space to the manifold at the point $\rho(y)$ and let $u_{1},
u_{2},..., u_{N-r}$ be the  components of $\clu(y):= (y-\rho(y))$ with respect to the
basis $\{e_{i}(y): i=1,\ldots,(N-r)\}$. 
We introduce a local coordinate system for the manifold $\cln_{\epsilon}P$ as follows:
$$G:\cln_{\epsilon}P\stackrel{F^{-1}}\rightarrow P\times
B_{\epsilon}^{(N-r)}(0)\stackrel{\xi\times {\rm id}}{\rightarrow}\IR^N~
({\xi~{\rm is~ a~ coordinate~ map ~for}~ P}):$$
$$y\rightarrow (\rho(y),\clu(y))\rightarrow( y_1,...y_r,u_1,...u_{N-r}  ).$$

Clearly, 
\be 
\label{coordinate_system} 
\lgl\lgl \frac{\partial}{\partial u_i}, \frac{\partial}{\partial u_j} \rgl\rgl = 0\ee   for all $i,j,~i=1,...,N-r, {\rm and} ~j=1,...,r$.

%
%
 Consider $\phi, \psi\in \widetilde{\cld}$ of the form $\phi=\xi\circ \rho$ and
$\psi=\beta \circ\clu$, where $\xi,\beta$ are smooth functions. Then
$\Phi(\phi)(y)=\eta(\xi)(\rho(y))$, $\Phi(\psi)=\psi\ot 1$, which gives by ( \ref{coordinate_system})  $$\lgl\lgl d\phi,d\psi\rgl\rgl=0,~~ 
\lgl\lgl d\Phi_{(1)}(d\phi),d\Phi_{(1)}(d\psi)\rgl\rgl=\Phi(\lgl \lgl d\phi, d \psi \rgl \rgl)=0.$$
 
 A general element of $\tilde{\cld} $ is a $\IC$-linear combination of the product functions of the form $(\xi \circ \rho) (\beta \circ \clu)$, 
 $\xi \in C^\infty(P),~\beta \in C^\infty(B^{N-r}_\epsilon(0))$. Hence it is enough to verify the following:\begin{displaymath} \lgl \lgl d\Phi_{(1)}(df_1), d\Phi_{(1)}(df_2) \rgl \rgl =\Phi(\lgl \lgl df_1,df_2 \rgl \rgl ) \
 for \ f_i=(\xi_i \circ \rho)(\beta_i \circ \clu).\end{displaymath} But we have $\Phi(f_i)(y)=\eta(\xi_i)(\rho(y))\beta_i(\clu(y))$ and 
 $\lgl \lgl d(\xi_i \circ \rho), d(\beta_j \circ \clu) \rgl \rgl =0$ for $i,j=1,2.$ Using this observation, Leibniz rule and the fact that $\eta$ is inner product preserving, 
 we complete the 
  proof of {\it Step 3}. \vspace{2mm}\\
{\it Step 4:}\\
Finally, it is easy to observe from the construction of $\Phi$ that it satisfies the condition regarding commutativity of the partial derivatives as in 
 the hypothesis of Corollary \ref{split_main}. Hence we conclude from that corollary that $\clq$ must be commutative. \qed

 Combining the above theorem with the techniques developed by Bhowmick and
Goswami in
\cite{qiso_comp} and by Joardar-Goswami in \cite{gos_joardar}, one can prove the following: 
\bcrlre 
The quantum isometry group of  a noncommutative manifold obtained by
cocycle twisting of a
classical, connected, compact Riemannian  manifold $M$ in the sense of \cite{rieffel}
is  a similar cocycle twisted version of $C(ISO(M))$, where $ISO(M)$ denotes the isometry
group of
$M$.
\ecrlre

 {\bf Acknowledgement:} 
 An initial version of this paper had Biswarup Das as a co-author. Later on the paper had undergone several rounds of revisions and corrections, in which Biswarup was not involved. 
 In his opinion, 
   his contribution to the present (final) draft does not remain significant enough to be a co-author and we have accepted his request to drop his name from the list of authors. However, 
    we would like to thank him for useful discussion in the initial stage of this work, particularly concerning Step 3 of the proof of Theorem \ref{main_thm}.

 We also thank professors Pavel Etingof, Chelsea Walton and Marc Rieffel and  Jyotishman Bhowmick and Subhrajit Bhattacharjee for encouragement, comments and discussion.

\end{document}